\documentclass{article}
\usepackage[utf8]{inputenc}
\usepackage{authblk}
\usepackage{amssymb}
\usepackage{amsthm}
\usepackage{amsmath}
\usepackage{amsxtra}
\usepackage{amsfonts}
\usepackage{mathtools}
\usepackage{bm}
\usepackage{siunitx}
\usepackage{tikz}
\usepackage[top=1in,bottom=1in,right=1in,left=1in]{geometry}
\usepackage{tkz-euclide}
\usepackage{listings}
\usepackage{amsthm}
\usepackage{algorithm}
\usepackage{pstricks-add}
\usepackage{colonequals}
\usepackage{graphicx} 
\usepackage{booktabs}
\usepackage{fancyhdr}
\usepackage{algorithm}
\usepackage{algpseudocode}
\usepackage{enumitem}
\usepackage{subcaption}
\usepackage{diagbox}
\usepackage[all]{xy}
\usepackage[normalem]{ulem}
\usepackage{amsmath}
\usepackage{makecell}

\DeclareMathOperator*{\argmin}{arg\,min}
\usepackage{color}

\title{Weighted inhomogeneous regularization for inverse problems with indirect and incomplete measurement data}
\author[1]{Bosu Choi\thanks{bchoi12@gsu.edu}}
\author[2]{Jihun Han\thanks{jihun.han@dartmouth.edu}}
\author[2]{Yoonsang Lee\thanks{yoonsang.lee@dartmouth.edu}}
\affil[1]{Department of Mathematics and Statistics, Georgia State University}
\affil[2]{Department of Mathematics, Dartmouth College}
\date{}

\begin{document}


\maketitle
\begin{abstract}
Regularization is a critical technique for ensuring well-posedness in solving inverse problems with incomplete measurement data. Traditionally, the regularization term is designed based on prior knowledge of the unknown signal’s characteristics, such as sparsity or smoothness. Inhomogeneous regularization, which incorporates a spatially varying exponent p in the standard $\ell_p$-norm-based framework, has been used to recover signals with spatially varying features. This study introduces weighted inhomogeneous regularization, an extension of the standard approach incorporating a novel exponent design and spatially varying weights. The proposed exponent design mitigates misclassification when distinct characteristics are spatially close, while the weights address challenges in recovering regions with small-scale features that are inadequately captured by traditional $\ell_p$-norm regularization. Numerical experiments, including synthetic image reconstruction and the recovery of sea ice data from incomplete wave measurements, demonstrate the effectiveness of the proposed method.
\end{abstract}

\section{Introduction}\label{sec:intro}
Inverse problems, which aim to recover an unknown signal from partial and noisy measurements, are ubiquitous in science and engineering. Examples include reconstructing objects using synthetic aperture radar (SAR) for remote sensing applications \cite{SARbook} and estimating sea ice thickness—an important geophysical parameter—from electromagnetic field data \cite{seaicethickness}. The CryoSat-2 satellite from the European Space Agency (ESA), for instance, measures the height difference between floating ice and surrounding water to infer polar ice properties \cite{cryosat2}.

A key challenge in inverse problems is the inherent lack of information needed for a unique reconstruction of the unknown signal. For instance, CryoSat-2 calculates sea ice thickness based on mean freeboard measurements over 25 km grids, effectively yielding a coarse spatial resolution of 125 km \cite{cryosat2resolution}. This limits the recovery of fine-scale features like cracks in sea ice. While NASA’s ICESat-2 satellite improves spatial resolution to a few meters, its coverage is limited to a small fraction of global ice floes \cite{ICESAT2}. Moreover, these measurements are indirect; sea ice thickness must be inferred through a nonlinear relationship between the measured height differences and the desired parameter.

To address such challenges, inverse problems are often formulated with additional constraints to ensure well-posedness, resulting in regularized optimization problems. For example, $\ell_2$ (Tikhonov) regularization is widely used for signals with smooth variations due to its computational simplicity. When the signal is sparse, $\ell_1$-based regularization is preferred, leveraging sparsity to improve reconstruction \cite{CS2006, CS2013}. Total Variation (TV) regularization extends this idea to recover signals with edges by enforcing sparsity in their gradients \cite{TV1992}.

Inhomogeneous regularization \cite{inhomoregv1} introduces a spatially varying exponent $p$ in $\ell_p$-norm constraints, enabling the recovery of signals with diverse features such as sparsity, smoothness, and oscillations. Unlike weighted $\ell_1$ regularization, which rescales the level set geometry, inhomogeneous regularization adapts the geometry itself by varying $p$ across the signal domain. While this approach has shown promise in recovering signals with mixed features, its performance degrades when features overlap spatially, causing misclassification. For instance, when oscillatory and discontinuous regions are adjacent, statistical patch-wise classification often fails, as observed in sea ice thickness recovery where edges, rapid variations, and smooth regions frequently coexist \cite{seaicemodeling, seaice}.

This study proposes weighted inhomogeneous regularization to overcome these limitations. Our approach combines new strategies for exponent and weight design, leveraging directional derivatives of smoothed reconstructions to refine classification and improve recovery. The directional derivatives clarify transitions between edges and oscillatory regions, while weights enhance the recovery of small-scale features by rescaling constraint geometries to incorporate additional prior knowledge. Furthermore, we restrict the range of $p$ to $[1, 2]$ to maintain convexity and solve the optimization problem using an adapted Alternating Direction Method of Multipliers (ADMM) algorithm \cite{ADMM}.

The paper is organized as follows: Section \ref{sec:standard_inhomo} reviews the inverse problem framework and the standard inhomogeneous regularization method \cite{inhomoregv1}. Section \ref{sec:newmethod} introduces the proposed weighted inhomogeneous regularization, detailing the exponent and weight design strategies and the ADMM-based solver. Section \ref{sec:numerics} demonstrates the method’s effectiveness through synthetic and real-world sea ice recovery experiments. Finally, Section \ref{sec:conclusion} discusses limitations and future research directions.

\section{Inverse problem and standard inhomogeneous regularization}\label{sec:standard_inhomo}

In this section, we examine the inhomogeneous regularization \cite{inhomoregv1} aimed at reconstructing multifeatured signals. We begin by outlining the general mathematical framework for a relevant inverse problem, followed by a discussion of the inhomogeneous regularization, addressing its rationale and associated challenges. The details are intentionally kept brief to allow for a clear understanding of the main concept behind the proposed method in the next section.

We consider the inverse problem expressed as follows:
\begin{equation}
	\bm{d} = G\bm{u}^* + \bm{\epsilon},
	\label{eq:problem}
\end{equation}
which involves retrieving an unknown signal $\bm{u}^* \in \mathbb{R}^n$ from an incomplete and indirect measurement $\bm{d}\in \mathbb{R}^m$ that is affected by measurement noise $\bm{\epsilon} \in \mathbb{R}^m$. Assuming that $G$, the measurement operator, is linear, it can be represented as a matrix in $\mathbb{R}^{m\times n}$. In scenarios where the measurement is incomplete, the size of the measurements, $m$, is generally less than the dimension $n$ of the unknown signal, making it challenging to obtain a unique solution. Additionally, the indirect measurement does not capture the direct elements of the signal. Therefore, any prior feature identification, such as an approximated shape of the target signal, is not available. An example of such an incomplete and indirect measurement is the subsampled Fourier measurement often encountered in remote sensing applications like synthetic aperture radar (SAR) \cite{SAR}.

Regularization introduces extra constraints on the target signal to enhance the well-posedness of the inverse problem. After selecting a specific regularization, $\mathcal{R}(\bm{u})$, the solution can be found by solving the regularized inverse problem within a restricted solution space, as represented by the following formulation:
\begin{equation}
	\bm{u}_{\lambda} = \argmin_{\bm{u} \in \mathbb{R}^n} \|\bm{d} - G\bm{u}\|_2^2 + \lambda \mathcal{R}(\bm{u}).
	\label{eq:opt_reg}
\end{equation} 
Here, $\lambda>0$ serves as a regularization parameter that balances the fidelity and regularization terms. The choice of regularization should be appropriate according to the characteristics of the target. For instance, Tikhonov or $\ell_2$ regularization, $\mathcal{R}(\bm{u})=\|\bm{u}\|_2$, is commonly used for signals with smooth variations, while $\ell_1$ norm regularization is favored for recovering sparse targets within a certain basis representation. Another type is Total Variation (TV) regularization, defined as $\mathcal{R}({\bm{u}})=\|D\bm{u}\|_1$, where $D$ denotes a discrete differential operator. By enforcing sparsity in the gradient domain, TV regularization can effectively recover targets characterized by edges or interfaces.

When signal solution exhibits multiple features without a predominant characteristic, traditional regularization methods utilizing homogeneous $\ell_p$ and TV norms may struggle to accurately recover it. While weighted homogeneous regularization offers some flexibility by allowing component-wise rescaling of constraints, it still faces limitations. For instance, weighted TV regularization can only recover a restricted range of piecewise continuous solutions. In cases where the solution contains both discontinuities and high-amplitude oscillations, weighted TV regularization does not effectively reconstruct the oscillatory sections. To address such mixed-type solutions, inhomogeneous regularization has been proposed in \cite{inhomo_denoising, inhomo_denoising2010, inhomoregv1}, which adapts the exponents of $\ell_p$ regularization according to local features. While the works in \cite{inhomo_denoising, inhomo_denoising2010} focus on denoising issues with direct measurements, \cite{inhomoregv1} investigates inverse problems with indirect measurements. Although it is feasible to estimate solution features from direct measurements, achieving the same from indirect measurements is more complex. Hence, in \cite{inhomoregv1}, the solution features are derived from the statistics of reconstruction samples. 

This approach leverages patch-wise statistical information and the spatial relationships between neighboring patches to identify local features. By dividing an image into smaller patches, we can analyze the statistical properties of these localized regions, reducing the uncertainty associated with pixel-wise statistics, which often suffer from noise and variability. This patch-wise approach enables more stable feature classification and exponent design. However, challenges arise when distinct features are closely spaced, potentially leading to confusion. To address this limitation, we integrate both patch-wise and pixel-wise approaches to improve the accuracy of feature classification, regardless of the distance between features. For a detailed explanation of our combined approach, please refer to Section \ref{sec:newmethod}, where we discuss concepts and details of standard inhomogeneous regularization from \cite{inhomoregv1}.

The formulation of inhomogeneous regularization is given by:
\begin{equation}
	\mathcal{R}(\bm{u}; \bm{p}) := \sum\limits_{j=1}^n |(D\bm{u})_j|^{p_j},
\end{equation}
where the spatially varying exponent $p_j$ is determined based on the characteristics of local features, and $(D\bm{u})_j $ denotes the $j$-th element of $ D\bm{u} $. In a one-dimensional signal, $ |(D\bm{u})_j| $ represents the magnitude of the discrete difference of $ \bm{u} $ at position $ j $, while in a two-dimensional image, $ |(D\bm{u})_j| $ can represent either isotropic or anisotropic total variation. 

To characterize the values of $ {p_j}$ for ${j=1, \cdots, n} $, \cite{inhomoregv1} employs patch-wise statistics to identify features and assign appropriate $ p_j $ values within each patch. The rationale for assigning exponents in a patch-wise manner, along with the estimation of patch-wise features, is that effective regularization occurs when a cluster of sufficiently large areas is assigned a constant exponent corresponding to the feature present. For example, if a discontinuity or edge occurs in a very small region, simply regularizing a few pixels precisely where the small feature resides with an exponent $ p_j=1 $ may not yield stable edge recovery. In contrast, recovery can be achieved by allocating an exponent of 1 to both the edge location and its surrounding area, thereby forming a sufficiently expansive exponent cluster.

As a result, we can express the inhomogeneous regularization incorporating patch-wise feature classification and exponent assignment from \cite{inhomoregv1} as follows:
\begin{equation}
	\mathcal{R}(\bm{u}; \bm{p}) = \sum\limits_{i=1}^M \sum\limits_{j \in I_i} |(D\bm{u})_j|^{p_j},
\end{equation}
where $ M $ denotes the number of patches, and $ I_i $ is the index set of pixels within the $ i $-th patch. The sets $ I_i $ are non-overlapping and their union, $ \{ I_i\}_{i=1}^M $, encompasses all pixels in the entire domain.

To estimate the local features of the solution that cannot be directly derived from indirect measurements, a reconstruction sample is employed. While each reconstruction may possess imperfections, statistical analysis allows us to create an approximate representation of the feature distribution. This sample is generated by solving multiple least-squares problems using classic homogeneous regularizations, expressed as follows:
\begin{equation}
\bm{u}_{\lambda_k; \bm{p}} = \argmin\limits_{\bm{u} \in \mathbb{R}^n} \| \bm{d}- G\bm{u} \|_2^2 + \lambda_k \mathcal{R}(\bm{u}; \bm{p}),
\label{eq:multRecons}
\end{equation}
for different choices of exponents, $p \in \{1, 2\}$, and regularization parameters $k = 1, 2, \cdots, C$. By adjusting the regularization parameter $\lambda_k$, we can compile a sample of size $C$ for each selection of $\bm{p}$. To derive the characteristics of the unknown signal, we initially use uniform values of 1 and 2 for the elements of $\bm{p}$, resulting in multiple reconstructions labeled as $\bm{u}_{\lambda_k; \bm{1}}$ and $\bm{u}_{\lambda_k; \bm{2}}$, respectively.

To obtain patch-wise statistics for feature identification, the standard inhomogeneous regularization utilizes the average gradient of sample reconstructions, defined as follows:
\begin{equation}
	\bm{g}_{\bm{p}} := \frac{1}{C} \sum_{k=1}^C D\bm{u}_{\lambda_k; \bm{p}},
	\label{eq:aveGrad}
\end{equation}
where $D$ is a discrete differential operator. Specifically, the variances of the average gradient components within an index set, $I_i$, are employed for feature classification, as they demonstrate greater stability and accuracy in distinguishing features compared to other statistical measures. In \cite{inhomoregv1}, features are categorized into three types: `discontinuity', `smoothness', and `oscillation'. Following the classification of features, the exponents $p_j$ are assigned based on the identified features. To create a sufficiently sized exponent cluster, all exponents $p_j$ corresponding to each $I_i$ are uniformly set within the range of [1,2]. If a patch is classified as `discontinuity,' the exponent cluster is assigned a value of 1 to encourage sparsity in the gradients. For a patch that shows no variation, i.e., it is flat, it is classified as `smoothness', and the exponent of 1 is also assigned. In other cases, exponents between 1 and 2 are allocated based on the levels of oscillation. As the level of oscillation increases, an exponent closer to 2 is more suitable for accurately capturing the oscillatory behavior. Conversely, when the signal is slightly changing or nearly flat, an exponent closer to 1 is more effective.

\begin{figure}[tb!]
	\makebox[\linewidth]{
		\includegraphics[width=\textwidth]{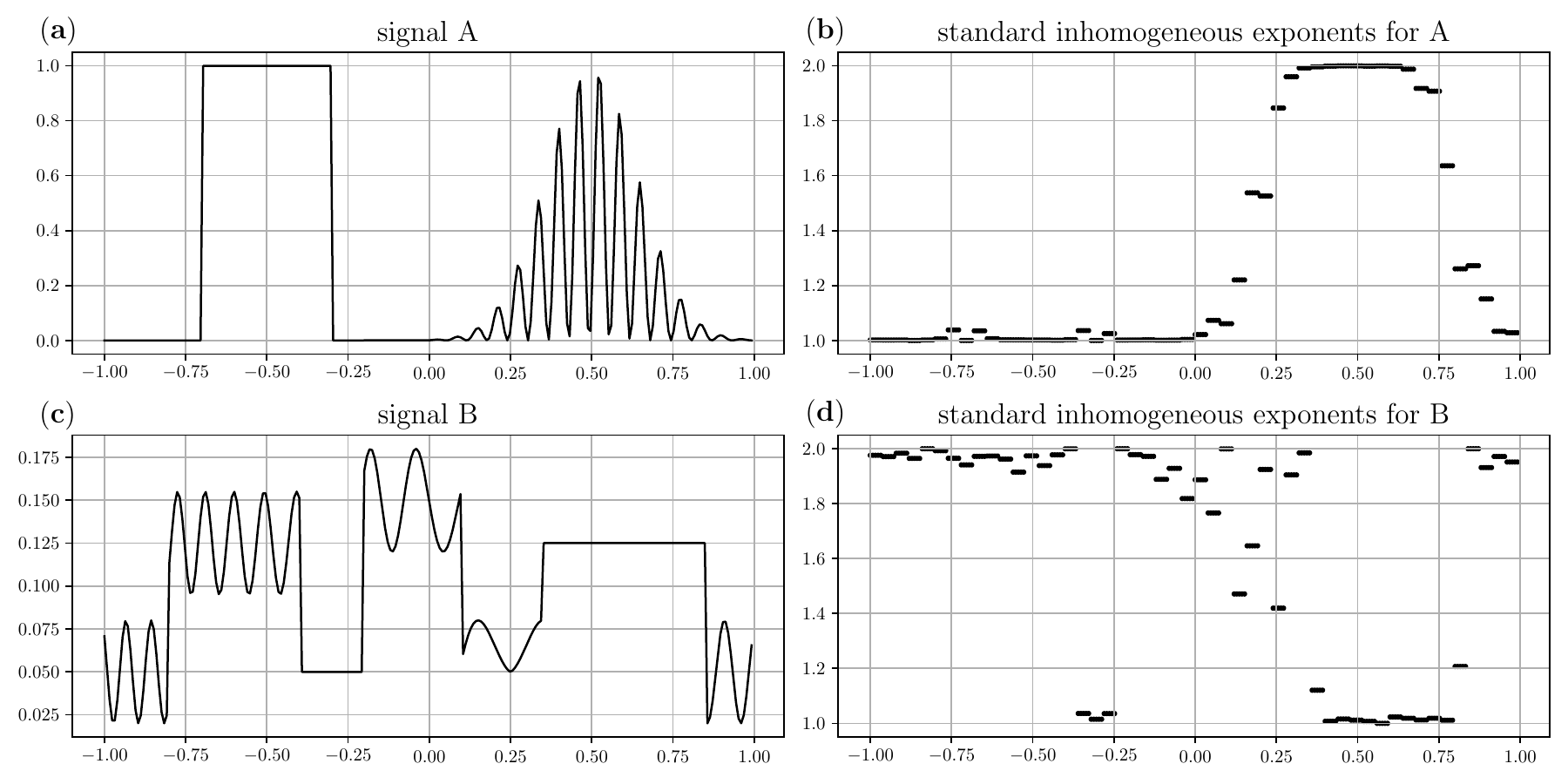}}
	\caption{True signals and exponent distributions from the standard inhomogeneous regularization
		\label{fig:standard_inhomo_exponent}}
\end{figure}

Inhomogeneous regularization based on patch-wise statistics successfully classifies features and assigns proper exponents if each feature is spatially separate with enough distance. Such an example is shown in Fig \ref{fig:standard_inhomo_exponent}(a), and the corresponding exponent design from the standard inhomogeneous regularization is shown in Fig \ref{fig:standard_inhomo_exponent}(b). We observe that 1 and values close to 1 are assigned around the discontinuities and smooth regions. On the other hand, values close to 2 are designated over the highly oscillatory region, and values gradually changing between 1 and 2 are set for the less oscillatory region. However, if different features are too close to each other, e.g., oscillation surrounding discontinuities, then the standard inhomogeneous regularization can fail to classify features correctly. Such an example is shown in Fig \ref{fig:standard_inhomo_exponent}(c), and the corresponding exponents from the standard inhomogeneous regularization are in Fig \ref{fig:standard_inhomo_exponent}(d). It is observed that the flat region takes exponents close to 1, and the oscillation takes exponents close to 2. However, edges adjacent to large amplitude oscillation are not well identified, which makes the exponents close to 2 assigned instead of 1. The assignment of improper exponent around edges happens because both center and neighborhood patches have large gradient variances, and therefore, patches with edges are misclassified as `oscillation.' 
Furthermore, the sea ice thickness from our motivation often contains such characteristics due to the dynamical movements of sea ice, including ice block breaks, pressure ridges, etc.

To resolve the misclassification by the standard inhomogeneous regularization from \cite{inhomoregv1}, we propose a new exponent design based on the combination of the patch-wise and pixel-wise feature estimation in Section \ref{sec:newmethod}. To enhance the efficacy of the updated inhomogeneous exponent for small feature recovery, we also introduce a weight design based on pixel-wise feature estimation. The weighted inhomogeneous regularization combining the new exponent and weight designs enables good recovery of multi-features of various sizes, which is demonstrated with synthetic and real sea ice images in Section \ref{sec:numerics}.

\section{Weighted inhomogeneous regularization}\label{sec:newmethod}

We present a new inhomogeneous exponent distribution that tackles important drawbacks of the standard inhomogeneous method covered in Section \ref{sec:standard_inhomo}. Our approach i) addresses feature misclassification and ii) integrates weighting techniques to improve the efficacy of inhomogeneous regularization. In order to improve classification accuracy, we utilize a new combination of patch-wise statistical analysis and pixel-wise estimation to rectify the misclassification of mixed features. To properly represent local feature fluctuations, the inhomogeneous exponents are adjusted based on these modifications. Solution recovery issues, however, could still exist, especially for features that are restricted to tiny areas and are not adequately highlighted by the exponents' limited influence. We address this by implementing a weighting technique that, when necessary, adds emphasis or penalization, improving the reconstruction of fine details and ensuring more robust feature recovery.

We emphasize the necessity of both patch-wise and pixel-wise statistics for accurate feature estimation, as more informative statistics vary depending on the feature type. Patch-wise statistics effectively identify smooth and oscillatory regions, as well as edges separate from oscillatory regions. Conversely, pixel-wise statistics are better at clearly identifying edges adjacent to oscillations, while the patch-wise approach may misclassify these edges as oscillations. Such misclassification arises because gradient variances of both center and neighboring patches become large for oscillatory regions lacking an edge and for edges next to oscillatory regions, which results in the patch-wise approach failing to adequately distinguish the two features. Therefore, these features are not effectively differentiated by patch-wise statistics. To overcome this challenge, we combine patch-wise and pixel-wise statistics to enhance the overall classification accuracy of various features.

Once features are accurately identified, we assign a uniform exponent value for each patch to ensure that each cluster of exponents is sufficiently large for effective regularization. This method may be less effective if two or more features are mixed within a patch. To mitigate the drawbacks of the patch-wise exponent design, we propose enhancing the information regarding feature variation within a patch through weighting, as the weight can indicate what should be emphasized or diminished within that patch. Consequently, our weight is formulated based on pixel-wise statistics rather than patch-wise statistics.

We formulate a general weighted inhomogeneous regularization, denoted as $\mathcal{R}(\bm{u}; \bm{\omega}, \bm{p})$, as follows:
\begin{equation}
	\mathcal{R}(\bm{u}; \bm{\omega}, \bm{p}) = \sum\limits_{i=1}^M \sum\limits_{j \in I_i} \omega_j |(D\bm{u})_j|^{p_j},
	\label{eq:weightedInhomoReg}
\end{equation}
where $\omega_j>0$ and $p_j \in [1,2]$ are the $j$-th elements of the weight and exponent distributions, $\bm{\omega}$ and $\bm{p} \in \mathbb{R}^n$, respectively. In Section \ref{sec:weightedinhomodesign}, we present our approach for determining $\bm{p}$ and $\bm{\omega}$ based on pixel-wise and patch-wise statistics. To more effectively differentiate edges and oscillations that are in close proximity, we incorporate directional derivatives of reconstructions convolved with a smoothing kernel. We address the feature misclassification from the patch-wise statistics using this pixel-wise estimation. Additionally, pixel-wise estimation is also employed to design a weight distribution, where its combination with a small cluster of inhomogeneous exponents aids in enhancing the reconstruction of small mixed features. In Algorithm \ref{alg:newinhomoreg}, we detail the complete procedure for designing the weighted inhomogeneous regularization and solving the regularized inverse problem. In Section \ref{sec:solver4weightedinhomodesign}, we introduce an optimization solver using ADMM to tackle a least-squares problem with weighted inhomogeneous regularization. The solver iteratively and alternately solves suboptimizations, decoupling the fidelity and regularization terms for computational efficiency.

\subsection{Weight and exponent design based on pixel-wise and patch-wise estimation} \label{sec:weightedinhomodesign}

Smooth and oscillatory regions are defined by the strength of variation within each region, leading us to consider patch-wise gradient variance instead of pixel-wise estimation. By examining the gradient variances in the center and neighboring patches, we can classify the oscillatory region more reliably, as gradient variances would remain consistently substantial. In contrast, the smooth region exhibits small gradient variances in both the center and neighboring patches. Therefore, we employ thresholding based on the size of patch-wise statistics to achieve stable classification of oscillation and smoothness. 

The patch-wise gradient variance is calculated as follows:
\begin{equation}
	var(\bm{g}_{\bm{p}}; I_i) = \frac{1}{|I_i|} \sum\limits_{j\in I_i} |g_{p,j}|^2 
	- \left( \frac{1}{|I_i|} \sum\limits_{j\in I_i} |g_{p,j}| \right)^2,
	\label{eqn:var}
\end{equation}
where $g_{p,j}$ denotes the $j$-th element of $\bm{g}_{\bm{p}}$ from \eqref{eq:aveGrad}, and $|I_i|$ indicates the size of the index set $I_i$. The patch-wise variances are normalized using min-max normalization and are denoted as $\widetilde{var}(\bm{g}_{\bm{p}}; I_i)$.

Using the normalized variances from the center and its neighboring patches, each patch can be identified as either `discontinuity,' `oscillation,' or `smoothness.' For smooth patches, a low variance is expected, while patches exhibiting oscillation or an edge are expected to have high variance. To differentiate between oscillation and an edge, we also consider the variances of neighboring patches. For oscillatory regions, we anticipate both center and neighborhood patches to have high variances, while for edges isolated from oscillations, the center patch shows high variance with low variances in neighboring patches. This variance-based feature classification typically yields stable and accurate results, although it may fail if an edge is surrounded by oscillations.

In contrast to smoothness and oscillation, discontinuities or edges occur within a very narrow region. Significant gradients typically characterize discontinuities under discretization, leading us to prioritize pixel-wise estimation over patch-wise statistics. Nevertheless, our multifeatured solutions may also encompass large amplitude oscillations with significant gradients. Thus, gradients computed directly from sample reconstructions may not effectively differentiate discontinuity from large-amplitude oscillation. To improve separation, we focus on their primary distinction: discontinuities exhibit a one-time significant intensity change, while oscillatory regions show periodic intensity variations. Therefore, we apply convolution with a kernel to smooth the reconstruction gradients, enhancing the distinction between these two features.

Convolution with an appropriate kernel can mitigate periodic intensity changes in the oscillatory region, reducing the oscillation amplitude toward zero. As a result, gradients also become relatively smaller. The kernel smooths out the discontinuity, but not to the same extent as the oscillation, since the intensity variation occurs only once without periodicity. One side adjacent to a discontinuity displays a higher average intensity, while the other side has a lower average intensity. Consequently, even after smoothing or averaging with kernels, the elevation change at the discontinuity can still be captured. The gradients at the discontinuity following smoothing provide a good estimate of the elevation change from one side to the other.

To define the gradient of the smoothed average reconstruction, we first denote the average reconstruction from the multiple reconstructions in \eqref{eq:multRecons} as follows: 
\begin{equation}
	\bm{r}_{\bm{p}} := \frac{1}{C} \sum_{k=1}^C \bm{u}_{\lambda_k; \bm{p}} \in \mathbb{R}^n.
	\label{eq:avgRecon}
\end{equation}
We focus on the average reconstruction since direct feature extraction from the indirect measurement is not feasible. To accurately differentiate edges from large-amplitude oscillations, we smooth the average reconstruction, $\bm{r}_{\bm{p}}$, by convolving it with a discretized compactly supported kernel, $\bm{k}$, as follows: 
\begin{equation}
	(\bm{r}_{\bm{p}} \ast \bm{k})i := \sum_{j} {r}_{{p},j} k_{i-j+1},
\end{equation}
where $(\bm{r}_{\bm{p}} \ast \bm{k})_i$ denotes the $i$-th element of the convolution of $\bm{r}_{\bm{p}}$ and $\bm{k}$, ${r}_{{p},j}$ represents the $j$-th element of $\bm{r}_{\bm{p}}$, and $k_{i-j+1}$ is the $(i-j+1)$-th element of the kernel $\bm{k}$. As the support size of $\bm{k}$ increases, the variations in the oscillatory region are canceled out more effectively, resulting in a clearer distinction between edges and oscillations. However, excessively large support sizes could overly smooth out all features. If the kernel's support size aligns with the oscillation's period, the oscillation can be effectively smoothed. Additionally, the kernel's shape can also influence the outcome of the convolution. Key conditions for a kernel to effectively average variations include: i) it should be compactly supported, ii) be nonnegative, iii) decay radially, and iv) have gradients with opposite signs on either side of the origin. We will further explore the sensitivity of our proposed exponent design to the choice of kernel after introducing our method for identifying missed edges in the following paragraph.

To identify edges, we calculate the derivative of the smoothed reconstruction, $\tilde{\bm{r}}_{\bm{p}} := \bm{r}_{\bm{p}} \ast \bm{k}$, and designate edges when the normalized derivative exceeds a threshold, $\tau$. For one-dimensional signals, we define a jump indicator, $\mathcal{J} \in \mathbb{R}^n$, as follows:
\begin{equation}
	\mathcal{J} := \left|\widetilde{D}\tilde{\bm{r}}_{\bm{p}} \right|,
	\label{eq:jumpInd1D}
\end{equation}
where $D\tilde{\bm{r}}_{\bm{p}} \in \mathbb{R}^n$ represents the finite difference of $\tilde{\bm{r}}_{\bm{p}}$, and $\widetilde{D}\tilde{\bm{r}}_{\bm{p}}$ is the normalized version of $D\tilde{\bm{r}}_{\bm{p}}$ using min-max normalization. 

For two-dimensional images, we consider the four directional derivatives along the vectors $v_1 = (1,0)$, $v_2 = (0,1)$, $v_3 = (1/ \sqrt{2},1/ \sqrt{2})$, and $v_4 = (-1/\sqrt{2},1/\sqrt{2})$ to capture edges in all directions. Consequently, we redefine the jump indicator, $\mathcal{J}$, as 
\begin{equation}
	\mathcal{J}:=
	\max \left\{
	\left| \widetilde{D}_{v_1} \tilde{\bm{r}}_{\bm{p}} \right|,
	\left| \widetilde{D}_{v_2} \tilde{\bm{r}}_{\bm{p}} \right|,
	\left| \widetilde{D}_{v_3} \tilde{\bm{r}}_{\bm{p}} \right|,
	\left| \widetilde{D}_{v_4} \tilde{\bm{r}}_{\bm{p}} \right|
	\right\},
	\label{eq:jumpInd2D}
\end{equation}
where $\widetilde{D}_{v_i}\tilde{\bm{r}}_{\bm{p}}$ is the normalized directional derivative of $\tilde{\bm{r}}_{\bm{p}}$ along $v_i$. If $\mathcal{J}$ exceeds $\tau$ at index $j$, we conclude there is an edge at $j$ and reclassify the patch corresponding to $I_i$ that includes $j$ as an edge location. This indicates that we continue to utilize the patch-wise exponent assignment even while employing pixel-wise estimation to identify missing features, thus preserving the effectiveness of the regularization. The updated inhomogeneous exponent, denoted as $\hat{\bm{p}}$, is derived from the classification corrected from the standard inhomogeneous regularization.

\begin{figure}[t]
	\makebox[\linewidth]{
		\includegraphics[width=\textwidth]{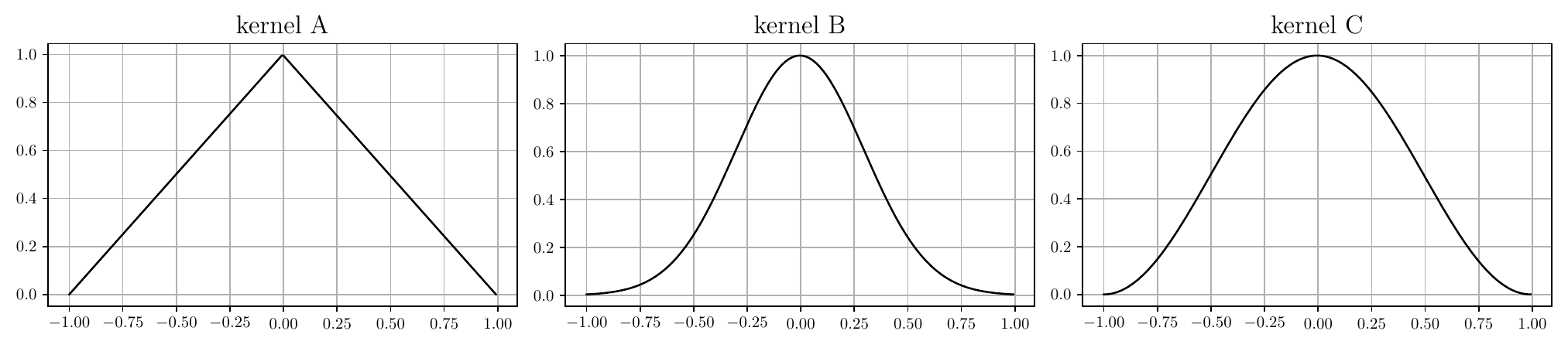}}
	\caption{Kernel A: triangular, Kernel B: Gaussian, Kernel C: shifted cosine
		\label{fig:kernels}}
\end{figure}

\begin{figure}[t]
	\makebox[\linewidth]{
		\includegraphics[width=\textwidth]{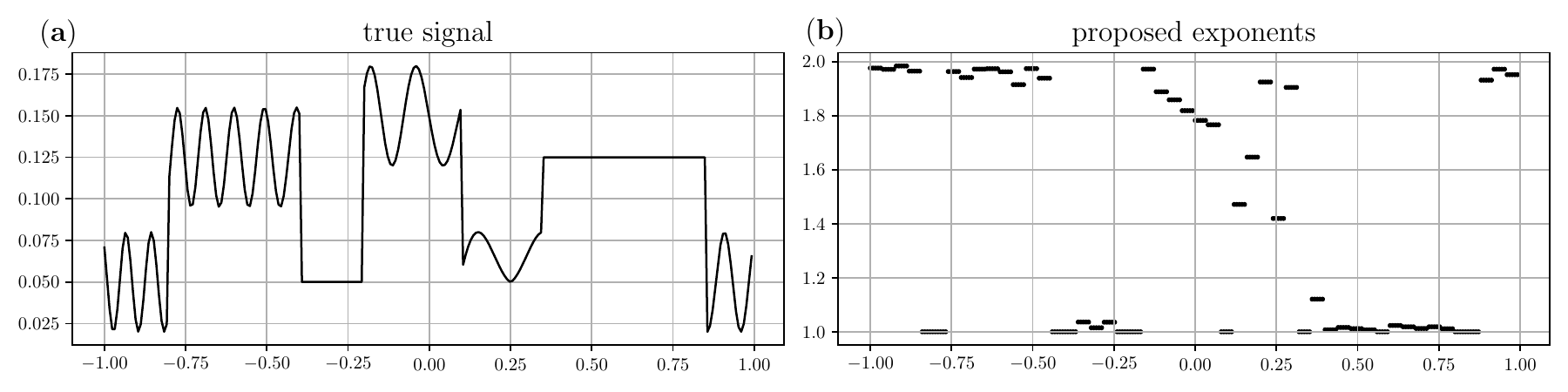}}
	\caption{True signal and exponent distributions from our proposed inhomogeneous regularization \label{fig:new_inhomo_exponent}}
\end{figure}

To inspect the sensitivity of the improved inhomogeneous exponent to different kernels with varying decay rates, we test three distinct kernels: triangular, Gaussian, and shifted cosine functions, as illustrated in Fig \ref{fig:kernels}. These kernels are scaled based on the selected vector support size, meaning they are adjusted to have smaller support when the vector support size is small. It is important to note that these kernels produce different forms of weighted averaging. Nevertheless, despite the varying averages produced by each kernel, we observe the same exponent design in Fig \ref{fig:new_inhomo_exponent}(b) resulting from the three kernels. This finding suggests that the resulting feature classification and corresponding exponent remain unaffected by the choice of kernel. Moreover, exponents are accurately assigned as 1 at all discontinuities, while the standard inhomogeneous exponents shown in Fig \ref{fig:standard_inhomo_exponent}(d) fail to capture discontinuities that are adjacent to oscillations. In addition to the discontinuities, smooth and oscillatory regions also correspond to appropriate exponents. This inhomogeneous exponent indicates that we have effectively updated the exponent design by integrating both pixel-wise and patch-wise statistics.

Once classification is complete, all exponents $p_j$ for $j$ belonging to $I_i$ are assigned uniformly based on the classification results. If a patch is categorized as 'discontinuity,' the exponent $p_j$ is set to 1 to encourage sparsity in the gradients. For oscillation and smoothness, the exponents are determined using the patch-wise average magnitude of gradients, defined as follows:
\begin{equation}
	avg(\bm{g}_{\bm{p}}; I_i) = \frac{1}{|I_i|} \sum\limits_{j\in I_i} |g_{p,j}|.
\end{equation}
We denote $\widetilde{avg}(\bm{g}_{\bm{p}}; I_i)$ as the normalized average using min-max normalization, and $p_j$ is defined as 
\begin{equation}
	p_j = 
	\begin{cases} 1 \hspace{15.3em} I_i \in {\rm discontinuity~ class},\\
			2 - \exp(-c \cdot \widetilde{avg}(\bm{g}_{\bm{1}}; I_i))
			\hspace{5em} I_i \in {\rm smooth~ class},\\
			2 - \exp(-c \cdot \widetilde{avg}(\bm{g}_{\bm{2}}; I_i))
			\hspace{5em} I_i \in {\rm oscillation~ class},
		\end{cases}
	\label{eq:standardInhomoExp}
\end{equation} 
where $c$ is a positive constant that ensures $p_j$ lies within the range $[1,2]$. If a target solution is highly oscillatory within a patch, resulting in a large average magnitude of the gradient, then a value close to $2$ is assigned to $p_j$. Conversely, if the solution is less oscillatory or smooth, a value close to $1$ is assigned.

In addition to the new inhomogeneous exponent, we design and integrate weights for the proposed weighted inhomogeneous regularizations. Although this new exponent represents the distribution of features more effectively, the quality of feature recovery can remain inadequate, especially when the target contains small features, as the cluster of corresponding exponents may be too small to be effective. In such cases, incorporating weights into the inhomogeneous regularization can enhance the efficacy of small exponent clusters. Furthermore, weights indicate which points should be emphasized or dampened if a patch with a constant exponent value encompasses multiple features. Therefore, to enhance the recovery quality of small features and mixed features within a patch, we determine weights pixel-wisely, in contrast to the patch-wise exponent distribution.

We design a weight distribution, $\bm{\omega}$, based on the derivative of the smoothed reconstruction, which accurately reflects the feature details. To support the recovery of various features, larger weights are anticipated for smooth regions, whereas relatively smaller weights are assigned for rapidly changing features such as edges and high-amplitude oscillations. Each weight is determined using a decreasing function $h(z)$ and the jump indicator $\mathcal{J}$ as follows:
\begin{equation}
	\omega_j := h(\mathcal{J}_j)
\end{equation}
where $\mathcal{J}_j$ represents the $j$-th element of $\mathcal{J}$. We utilize a linearly decreasing function $h(z) = az + b$ for weight design, which consistently reflects the intensity variations. The constants $a < 0$ and $b > 0$ are chosen to ensure that $h(z)$ is positive for $z \in [0,1]$. We keep $h(z)$ strictly greater than 0 since a weight of zero would nullify the regularization at the corresponding point. Furthermore, we maintain a relatively small range of weight values to achieve a balanced distribution that promotes all features.

Weighted homogeneous and inhomogeneous regularizations target different interests, leading to different strategies for weight design. Weighted $\ell_1$ and TV regularizations typically focus on recovering sparse and piecewise smooth functions. These functions generally consist of smooth regions that are constant or slowly varying. Conversely, our weighted inhomogeneous regularization accounts for solutions that include highly oscillatory regions as well as discontinuities and smoothness. While weighted $\ell_1$ and TV regularizations aim to identify sparse or edge locations rather than estimating overall variation, they penalize both constant and slowly varying regions significantly and promote sparsity or discontinuity using rapidly decaying functions, such as exponential or power functions. In contrast, our approach employs linearly decreasing functions and derivative estimation to set the weights for inhomogeneous regularization, reflecting all variations evenly. Additionally, the weight range for weighted homogeneous regularization tends to be broader than ours, leading to a stronger promotion of certain features. Consequently, weighted homogeneous regularization assigns a large weight to smooth areas while applying a minimal weight to sparsity and discontinuity, resulting in a dominant reconstruction of a specific feature. Our approach maintains a smaller weight range to avoid excessive promotion or penalization of any particular feature.

In Section \ref{sec:numerics}, we compare the performance of our proposed weighting strategy against another weighting method combined with multiple exponents. We find that our weighting strategy excels in recovering diverse types of features collectively. In Section \ref{subsec:synthetic2d}, our weighted regularization is applied to both unscaled and scaled features—larger and smaller features—to demonstrate that weighting enhances the regularization effect when a feature is small, resulting in a smaller exponent cluster. Thus, the proposed weighted inhomogeneous regularization proves superior in recovering targets with multiple features of varying sizes. In Section \ref{subsec:exp:seaice}, this is further validated by results showing that the proposed method effectively recovers real sea ice images containing complex features of different sizes.

\begin{algorithm}[t]
	\caption{Signal recovery with a weighted inhomogeneous regularization}\label{alg:newinhomoreg}
	\begin{algorithmic}
		\Require kernel support size $l$, threshold $\tau$, constants $a$ and $b$ 
		\State 1. Obtain multiple reconstructions for $k = 1, 2, \cdots, C$:
		$$ \bm{u}_{\lambda_k; \bm{p}} = \argmin\limits_{\bm{u} \in \mathbb{R}^n} \| \bm{d}- G\bm{u}\|_2^2 + \lambda_k \mathcal{R}(\bm{u}; \bm{p}),~ \bm{p}=\bm{1}, \bm{2} $$
		\State 2. Obtain the pre-classification using patch-wise statistics in \eqref{eqn:var} as in \cite{inhomoregv1}
		\State 3. Compute the average reconstruction, $\bm{r}_{\bm{p}}$ with $\bm{p}=\bm{1}$, as in \eqref{eq:avgRecon}
		\State 4. Compute the jump indicator, $\mathcal{J}$, as in \eqref{eq:jumpInd1D} for a 1D target and \eqref{eq:jumpInd2D} for a 2D target from the derivatives of the average reconstruction convolved with a kernel whose support size is $l$ 
		\State 5. Re-classify an index set $I_i$ to `discontinuity' if $I_i$ contains any index $j$ such that $\mathcal{J}_j  > \tau$
		\State 6. Design an exponent distribution $\hat{\bm{p}}$ based on the pre- and re-classifications by determining $p_j$ according to \eqref{eq:standardInhomoExp}
		\State 7. Compute each weight $\omega_j = h(\mathcal{J}_j) = a \cdot \mathcal{J}_j +b$
		\State 8. Recover the target by solving a least-squares problem with the weighted inhomogeneous regularization:
		$$ \bm{u}_{\lambda; \bm{\omega}, \hat{\bm{p}}} = \argmin\limits_{\bm{u} \in \mathbb{R}^n} \| \bm{d}- G\bm{u}\|_2^2 + \lambda \mathcal{R}(\bm{u}; \bm{\omega}, \hat{\bm{p}}) $$
	\end{algorithmic}
\end{algorithm}

We conclude this section by examining the pseudocode in Algorithm \ref{alg:newinhomoreg} for signal recovery using a weighted inhomogeneous regularization, $\mathcal{R}(\bm{u}; \bm{\omega}, \hat{\bm{p}})$. The algorithm requires several parameters: the kernel support size, $l$, for calculating the derivative of the smoothed reconstruction; a threshold, $\tau$, for reclassification; and constants, $a$ and $b$, which define the linear weight function. In Line 1, we gather multiple reconstruction samples, and in Line 2, we obtain the feature (pre)classification from \cite{inhomoregv1} using patch-wise statistics primarily to identify smoothness and oscillations. Lines 3-6 of Algorithm \ref{alg:newinhomoreg} correct any classification missed by the standard inhomogeneous regularization through pixel-wise estimation, focusing on identifying discontinuities and subsequently determining an exponent distribution $\hat{\bm{p}}$. In Line 7, we compute the weight that reflects the feature details. Finally, Line 8 addresses the least-squares problem with the weighted inhomogeneous regularization. The specific solver for the least-squares problem with the weighted inhomogeneous regularization in Line 8, along with the equivalent problem for homogeneous regularization that assigns uniform $p_j$ and $\omega_j$ in Line 1, is discussed in detail in Section \ref{sec:solver4weightedinhomodesign}.

\subsection{Solver for weighted inhomogeneous regularization} \label{sec:solver4weightedinhomodesign}
As the inhomogeneous distribution of regularization exponents is designed to have a  lower bound $1$ (i.e., $p_i \geq 1$, $\forall i$), a convex optimization method is appropriate for solving the problem. The alternating direction method of multipliers (ADMM; \cite{ADMM}) inherent to the augmented Lagrangian method decomposes the constraint optimization into a series of easier-to-handle unconstraint subproblems. In this paper, we employ ADMM to solve the proposed regularized optimization through a suitable problem reformulation, taking the computational advantages of decoupling the regularization term from the fidelity term and in parallel computing for multiple reconstructions. Also, inspired by the group Lasso regularization \cite{groupLasso}, we reduce the suboptimization associated with the regularization term separable, thereby averting the computational challenges caused by the coupling of inhomogeneous distributions of both weights and exponents in the proposed regularization.   

We reformulate the problem with the proposed regularization in \eqref{eq:weightedInhomoReg} into the constraint optimization as follows:
\begin{align}
	\text{minimize} ~~ \mathcal{Q}(\bm{u}) + \lambda \mathcal{R}(\bm{v}) \label{eq:ADMM Objectives}\\
	\text{subject to} ~~ \bm{F}\bm{u} - \bm{v} = \bm{0}, \label{eq:ADMM constraint}
\end{align}
where the objective functions $\mathcal{Q}$ and $\mathcal{R}$ are 
\begin{align}
	\mathcal{Q}(\bm{u}) &=  \| \bm{A}\bm{u}- \bm{y} \|^2_2, \nonumber \\
	\mathcal{R}(\bm{v}) &= \sum \limits_{i=1}^{n}g_i(\|\bm{v}_i\|_2), ~~~~ g_i(z) :=  \omega_i \alpha_i(z),~ \alpha_i(z) :=|z|^{p_i}.\label{eq:rerrange_reg}
\end{align}
We note that we introduce the auxiliary variable $\bm{v} =(\bm{v}_1, \bm{v}_2,\cdots, \bm{v}_l)\in \mathbb{R}^{l}$ with the constraint \eqref{eq:ADMM constraint}, in which $\bm{F} : \mathbb{R}^{n} \rightarrow \mathbb{R}^{l}$ is the discrete gradient operator. It is read as $\bm{v}_i =D {u}_i \in \mathbb{R},  i=1,2,\cdots, n$,  for the case of 1D signal (in case $l=n$),  and $\bm{v}_i = (D_x {u}_i, D_y {u}_i) \in \mathbb{R}^2,  i=1,2,\cdots, n$,  for the case of 2D image (in case $l=2n$). The augmented Lagrangian corresponding to the optimization \eqref{eq:ADMM Objectives} and \eqref{eq:ADMM constraint} is
\begin{equation}
	L_{\rho}(\bm{u}, \bm{v},  \bm{\eta}) = \mathcal{Q}(\bm{u}) + \lambda\mathcal{R}(\bm{v}) +  \bm{\eta}^{T}(\bm{F}\bm{u} - \bm{v}) + (\rho/2)\|\bm{F}\bm{u} - \bm{v}  \|_2^2,
\end{equation}
where $ \bm{\eta}\in \mathbb{R}^{l}$ is the Lagrange multiplier and $\rho>0$ is the penalty parameter.
ADMM alternately updates each variable $\bm{u}$, $\bm{v}$ and $ \bm{\eta}$ with the following series of iterations
\begin{equation}
	\bm{u}^{k+1} =\text{Prox}_{\mathcal{Q},\rho,\bm{F}}(\bm{v}^{k} -  \bm{\eta}^{k}):=\operatorname*{argmin}_{\bm{u}} \left(\mathcal{Q}(\bm{u})+ \frac{\rho}{2}\|\bm{F}\bm{u} -(\bm{v}^{k} -  \bm{\eta}^{k}) \|_2^2\right), \label{eq:u-proximal}
\end{equation}
\begin{equation}
	\bm{v}^{k+1} =\text{Prox}_{\mathcal{R},\frac{\rho}{\lambda},\bm{I}}(\bm{F}\bm{u}^{k+1} +  \bm{\eta}^{k}):=\operatorname*{argmin}_{\bm{v}} \left( \mathcal{R}(\bm{v}) + \frac{\rho/\lambda}{2}\| \bm{I}\bm{v} - (\bm{F}\bm{u}^{k+1}  +  \bm{\eta}^{k}) \|_2^2\right),\label{eq:v-proximal}
\end{equation}
\begin{equation}
	\bm{\eta}^{k+1} =  \bm{\eta}^{k}+ \bm{F}\bm{u}^{k+1} - \bm{v}^{k+1}.
\end{equation}
The suboptimizations,  \eqref{eq:u-proximal} and \eqref{eq:v-proximal} are known as proximal operators, which,  in general,  have either closed-form or implicit form depending on the objective terms.  The $\bm{u}$-optimization \eqref{eq:u-proximal},  as a quadratic objective form,  is known to have the  closed-form
\begin{equation}
	\text{Prox}_{\mathcal{Q},\rho,\bm{F}}(\bm{x})= \left(\bm{A}^T\bm{A}+\rho\bm{F}^{T}\bm{F})^{-1}(\bm{A}^T\bm{b}+\rho\bm{F}^{T}\bm{x}\right), ~~\bm{x} \in \mathbb{R}^{l}.
\end{equation}
For the proximal operator in $\bm{v}$-optimization \eqref{eq:v-proximal},  the objective function $\mathcal{R}(\bm{v})$ is separable,  so is the optimization.  The proximal $\text{Prox}_{\mathcal{R},\rho,\bm{I}}(\bm{x})$ is
\begin{align}
	\text{Prox}_{\mathcal{R},\rho,\bm{I}}(\bm{x}) &= \operatorname*{argmin}_{\bm{v}}\left(\mathcal{R}(\bm{v})+\frac{\rho}{2}\|\bm{v} - \bm{x} \|_2^2\right) \\
	&= \operatorname*{argmin}_{\bm{v}}\left(\sum \limits_{i=1}^{n}g_i(\|\bm{v}_i\|_2)+\frac{\rho}{2}\sum \limits_{i=1}^{n}\|\bm{v}_i - \bm{x}_i \|_2^2\right) \nonumber  \\
	&= \left(\cdots, \operatorname*{argmin}_{\bm{v}_i} \left( g_i(\|\bm{v}_i\|_2) + \frac{\rho}{2}\| \bm{v}_i - \bm{x}_i \|_2^2 \right), \cdots \right) \in \mathbb{R}^{l},  \nonumber
\end{align} 
and each component of the proximal is 
\begin{align}
	\left(\text{Prox}_{\mathcal{R},\rho,\bm{I}}(\bm{x})\right)_i
	&= \operatorname*{argmin}_{\bm{v}_i} \left( g_i(\|\bm{v}_i\|_2) + \frac{\rho}{2}\| \bm{v}_i - \bm{x}_i\|_2^2 \right) \nonumber  \\
	&= \operatorname*{argmin}_{\bm{v}_i} \left(  \omega_i \alpha_i (\|\bm{v}_i\|_2) + \frac{\rho}{2}\| \bm{v}_i - \bm{x}_i \|_2^2 \right) \nonumber  \\
	&= \operatorname*{argmin}_{\bm{v}_i} \left(  \alpha_i(\|\bm{v}_i\|_2) + \frac{\rho/  \omega_i}{2}\| \bm{v}_i - \bm{x}_i \|_2^2 \right) \nonumber \\
	&=\text{Prox}_{ \alpha_i(\|\cdot \|_2), \frac{\rho}{ \omega_i}, \bm{I}}(\bm{x}_i)\nonumber \\
	&= \text{Prox}_{ \alpha_i, \frac{\rho}{ \omega_i}, \bm{I}}(\|\bm{x}_i\|_2)\frac{\bm{x}_i}{\|\bm{x}_i\|_2},  ~~ i=1,2,\cdots,n, \label{eq:prox_norm_composition} 
\end{align}
in which the last equality holds by the norm composition rule of the proximal operator.  The proximal operator,  $\text{Prox}_{ \alpha_i, \frac{\rho}{ \omega_i}, \bm{I}}(q)$ in \eqref{eq:prox_norm_composition} is easily derived as 
\begin{equation}
	\text{Prox}_{ \alpha_i, \frac{\rho}{ \omega_i}, \bm{I}}(q) =
	\begin{cases}
		\mathcal{S}_{ \omega_i/\rho}(q) & ~\text{if } p_i =1,
		\\
		\text{The zero of }~ T(x) :=  \text{sgn}(x)p_i|x|^{p_i-1} + \frac{\rho}{ \omega_i} (x - q) & ~\text{if } p_i > 1,
	\end{cases}
\end{equation}
where $\mathcal{S}_{\kappa}$ is the soft thresholding or shrinkage operator defined as $\mathcal{S}_{\kappa}(q):=\left(1-\frac{\kappa}{q} \right)_{+}q$.  We note that $T(x)$ is an increasing function, and we find the unique root of $T(x)$ with Chandrupatla's method \cite{chandrupatla}.

\section{Numerical Experiments}\label{sec:numerics}


We showcase the effectiveness of our proposed weighted inhomogeneous regularization in reconstructing solutions with mixed features of diverse shapes and sizes through extensive numerical experiments on both synthetic and real sea ice images. Our method refines the inhomogeneous exponent distribution and produces a weight distribution that accurately captures the overall intensity variations. To validate the robustness of our approach, we compare it against classic homogeneous regularization and standard inhomogeneous regularization, both with and without weights, in recovering features of varying sizes.

For this comparison, we conduct numerical experiments on synthetic images and their scaled-down counterparts, featuring mixed and scaled features. Additionally, we test our approach on a real sea ice image, which presents more complex patterns and diverse feature sizes. To evaluate the performance of each method, we calculate the relative $\ell_p$ reconstruction errors, defined as:
\begin{equation}
	\frac{\| \hat{\bm{u}} - \bm{u}^* \|_p}{\| \bm {\bm{u}}^* \|_p},
\end{equation}
where $ \hat{\bm{u}}$ is a recovery of true $\bm{u}^* $. Moreover, we compare reconstructions and pointwise errors to see if various local features are well recovered.

We focus on the recovery of target features from indirect and incomplete measurements. Specifically, we use partial Fourier measurements, which are commonly employed in applications such as MRI and SAR data. For synthetic image tests (Section \ref{subsec:synthetic2d}), we assume noiseless measurements by setting $\bm{\epsilon} = \bm{0}$ in \eqref{eq:problem} to highlight our method's ability to generate accurate exponent and weight distributions representing feature variations. For real sea ice image tests (Section \ref{subsec:exp:seaice}), we incorporate additive Gaussian noise into the measurements to demonstrate the robustness of our method in noisy scenarios.

To estimate features from limited and indirect data, we generate sample reconstructions by solving least-squares problems with homogeneous regularizations. We employ $p=1$ and $p=2$ for the regularizations, varying the parameter $\lambda$ over a logarithmically spaced range within a random subinterval of length $10^2$ to $10^4$ in $[10^{-3}, 10^{3}]$. The average reconstruction derived from these samples is relatively insensitive to the specific subinterval choice, effectively capturing the target solution's characteristics. In our experiments, we use a sample size of $C=50$ for 2D image recovery.

Parameters required for the weighted inhomogeneous regularization are detailed in Algorithm \ref{alg:newinhomoreg}, with specific values for each test image provided in the corresponding subsections. Each least-squares problem with regularization is solved using the ADMM method described in Section \ref{sec:solver4weightedinhomodesign}, run in parallel for multiple $\lambda$ values. For all experiments, the penalty parameter $\rho$ in ADMM is set to 1.

To evaluate the effectiveness of our weighting method, which is derived from the derivatives of smoothened average reconstructions, we also compare it against an alternative strategy based on variance-based joint sparsity (VBJS) introduced in \cite{VBJS}. The VBJS method is particularly effective at identifying edges in relatively smooth solutions, except near discontinuities. It calculates weights using the pointwise variance, $\hat{var}(\cdot)$, of multiple reconstructions as follows:
\begin{equation}
	\hat{var}( P \bm{u}_{\lambda; {p}, j}) =\frac{1}{C} \sum_{k=1}^C 
	\left((P \bm{u}_{\lambda_k; \bm{p}})_j \right)^2 -  \left( \frac{1}{C} \sum_{k=1}^C  (P \bm{u}_{\lambda_k; \bm{p}})_j \right)^2 , ~j =1, 2, \cdots, n, 
\end{equation}
where $P$ represents a sparsifying operator, such as the polynomial annihilation operator from \cite{archibald2005polynomial}, and $ (P \bm{u}_{\lambda_k; \bm{p}})_j$ is the $j$-th element of $P \bm{u}_{\lambda_k; \bm{p}} \in \mathbb{R}^n$.
The weights, $\hat{\bm{\omega}}$, are computed as:
\begin{equation}
	\hat{w}_j = \frac{1}{\hat{var}( P \bm{u}_{\lambda; {p}, j}) + \hat{\epsilon}},
	\label{eq:vbjs}
\end{equation}
where $\hat{\epsilon} > 0$ is a small constant. In \cite{VBJS}, $\hat{\epsilon}$ is set to $0.01$, ensuring that the weights are confined to the range $[0, 100]$. Points with lower variance are assigned larger weights, while those with higher variance (typically observed at discontinuities) receive smaller weights. These smaller weights reduce the penalization at discontinuities, thereby promoting the corresponding elements in the reconstruction.

\subsection{2D synthetic image recovery}\label{subsec:synthetic2d}

\begin{figure}[tb!]
	\makebox[\linewidth]{
		\includegraphics[width=\textwidth]{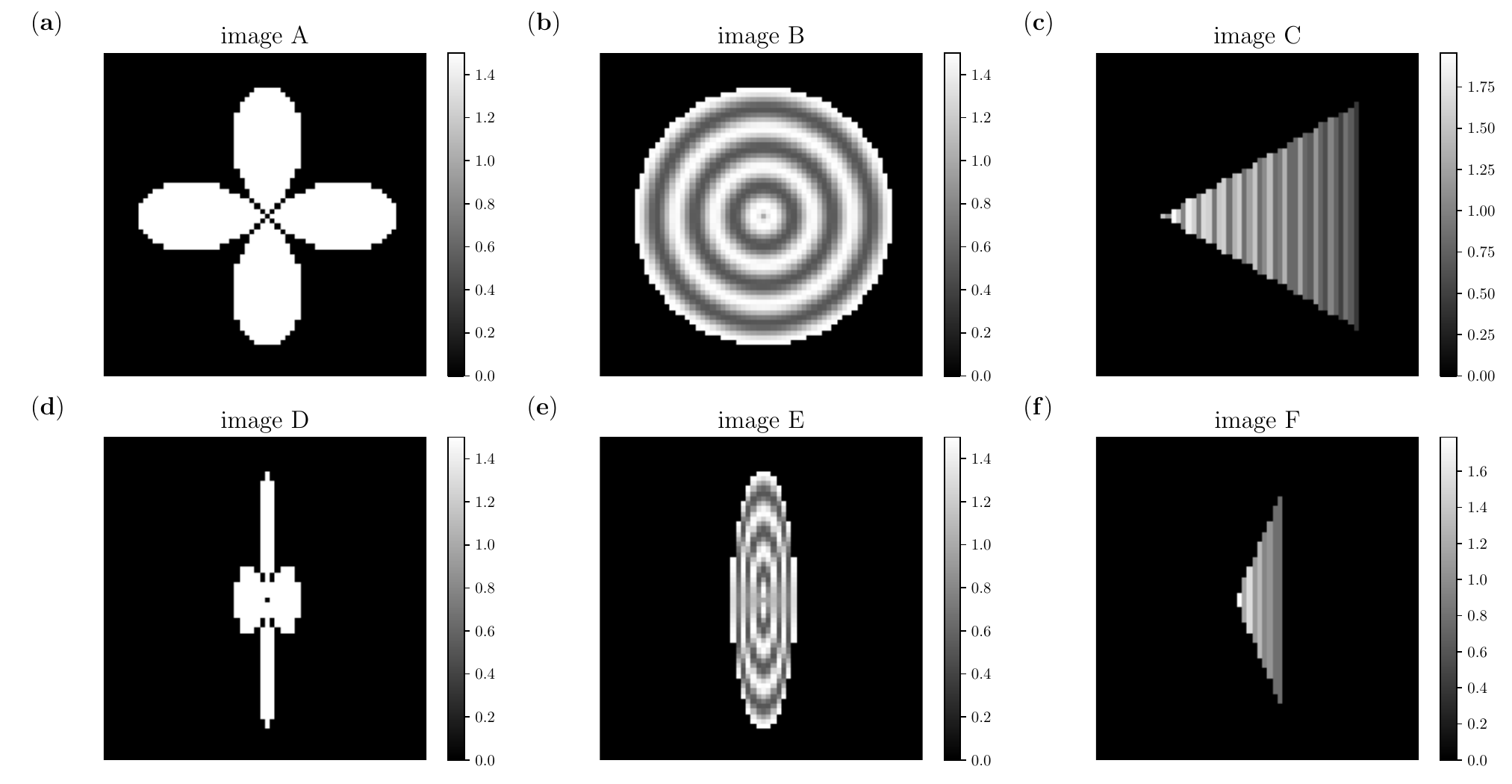}}
	\caption{2D synthetic images (A-C) and their scaled images (D-F) by a factor of 1/4 along $x$-axis
		\label{fig:true_2d_synthetic}}
\end{figure}

To evaluate the accuracy and effectiveness of our proposed exponent and weight design, we consider a variety of 2D synthetic images with mixed features of different shapes and sizes, as illustrated in Fig. \ref{fig:true_2d_synthetic}. Image A in Fig. \ref{fig:true_2d_synthetic}(a) consists of four leaves with curved edges, representing a piecewise constant image. Such images are typically well reconstructed using homogeneous $p=1$ regularization. Image B contains three distinct types of features: a circular edge, oscillations inside, and smooth areas outside. Notably, the proximity of the oscillations to the edge can lead to misclassification when using standard homogeneous regularization. Image C features a triangular ramp with internal oscillations, causing fluctuating edge heights. Unlike Image B, the oscillations in Image C propagate directionally rather than radially.
These three images exhibit a variety of feature shapes, such as edges along curves or lines, constant or varying edge heights, and oscillations with radial or directional variations. Our approach demonstrates the ability to accurately design exponents and weights that describe these complex features.

Additionally, we test scaled-down versions of these images, labeled as Images D-F in Fig. \ref{fig:true_2d_synthetic}, created by compressing Images A-C along the $x$-axis by a factor of 1/4. This scaling emphasizes the importance of combining inhomogeneous exponents and weights for recovering small features. Each pair of images (A and D, B and E, C and F) shares features of similar shapes but at different scales. Interestingly, the preferences for regularization design differ based on whether the features are scaled, as we demonstrate in subsequent sections.
For these tests, we discretize the images on a $64 \times 64$ uniform grid and use 39.1\% of the lowest frequency coefficients in the $y$-direction for measurements. This setup allows us to rigorously validate the adaptability and robustness of our approach.

\subsubsection{Recovery of images with large features}

We select Images A, B, and C as the initial test subjects to validate the capability of our proposed exponents in enhancing the standard inhomogeneous exponents by addressing misclassified features. Additionally, we demonstrate that while the proposed weights effectively capture the variations in the target features, they do not significantly improve the recovery quality for adequately large features. Furthermore, we analyze the impact of combining the proposed weights with various exponent designs on the overall recovery performance, comparing these results to those achieved using alternative weighting strategies.

\begin{figure}[htb!]
	\makebox[\linewidth]{
		\includegraphics[width=\textwidth]{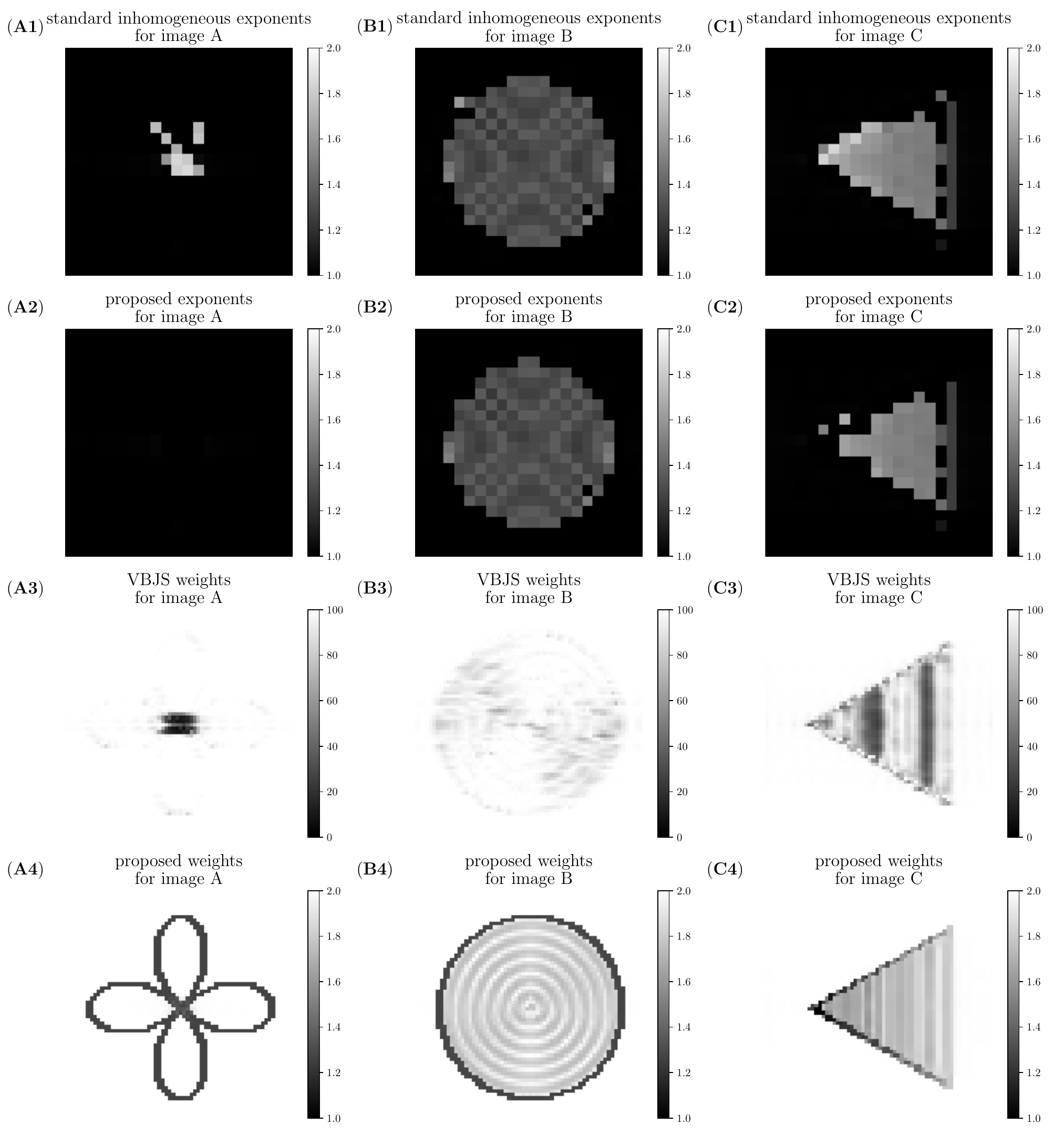}}
	\caption{Standard and proposed exponents, and VBJS and proposed weights for Images A-C in Fig \ref{fig:true_2d_synthetic}(a)-(c)
		\label{fig:exponent_weight_ABC}}
\end{figure}

\begin{figure}[tb!]
	\makebox[\linewidth]{
		\includegraphics[width=\textwidth]{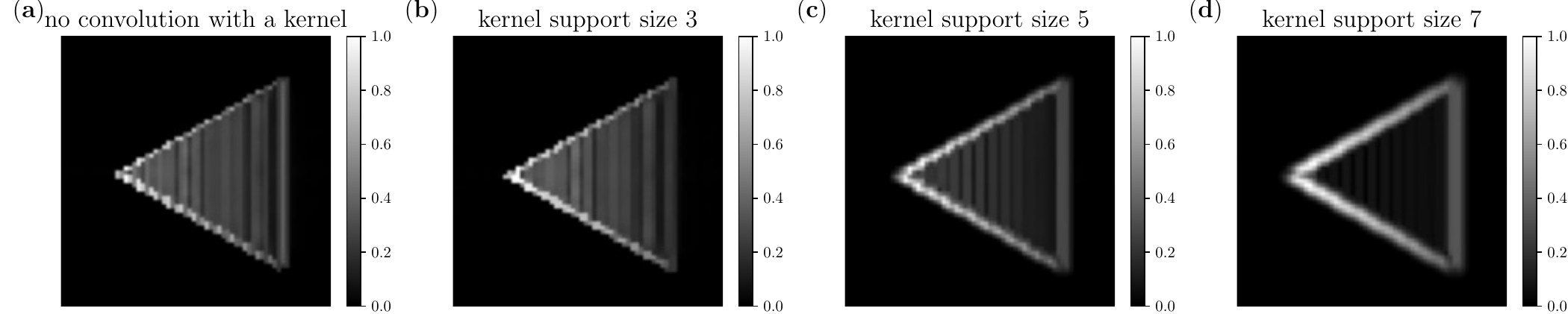}}
	\caption{Comparison of jump indicators of unsmoothed and smoothened average sample reconstructions for Image C in Fig \ref{fig:true_2d_synthetic} with kernel support size 3, 5, and 7
		\label{fig:jump_ind_comp}}
\end{figure}

When comparing the standard and proposed inhomogeneous exponents for Images A-C in Fig. \ref{fig:exponent_weight_ABC}, a recurring pattern emerges. The standard inhomogeneous exponents tend to misclassify portions of edges as oscillations, assigning relatively large exponent values instead of 1. In contrast, the proposed exponents correct these misclassifications and assign the appropriate value of 1 to the corresponding areas.
This misclassification by the standard inhomogeneous exponents typically occurs in regions where the patch-wise gradient variances are high both at the center and in the vicinity of edges. For instance, in Image A, the central area contains intersecting edges, leading to several patches with high gradient variances. These patches are misclassified as oscillations rather than discontinuities due to the statistical similarity between true oscillations and the observed patterns. As a result, the central region shows large exponent values, as depicted in Fig. \ref{fig:exponent_weight_ABC}(A1).

For Images B and C, which feature curved and straight edges adjacent to oscillatory regions, the patches covering both the edges and oscillations exhibit large gradient variances. This overlap leads to misclassification of the edges bordering the oscillatory regions. Consequently, certain parts of these edges are assigned large exponents, as observed in Fig. \ref{fig:exponent_weight_ABC}(B1) and (C1).
In contrast, the proposed exponents in Fig. \ref{fig:exponent_weight_ABC}(A2)-(C2) accurately identify edges that were previously misclassified using the patch-wise variance-based approach, correctly assigning an exponent of 1 to those regions. Additionally, we observe that the "inhomogeneous" exponents in Fig. \ref{fig:exponent_weight_ABC}(A2) are not uniformly 1. Instead, while patches containing edges have exponents exactly equal to 1, those covering smooth regions have exponents slightly greater than 1, reflecting their distinct characteristics.

To identify edges missed by the standard inhomogeneous regularization, the derivatives of the smoothened average sample reconstruction prove effective in better distinguishing edges from oscillations. By applying convolution with a kernel vector, the derivatives over oscillatory regions are flattened, while the derivatives at edges remain relatively pronounced. The smoothing effect depends on the kernel support size, with larger support sizes providing greater flattening of oscillations. However, as long as the kernel satisfies the key conditions (i)-(iv) outlined in Section \ref{sec:weightedinhomodesign} and the support size remains relatively small, the convolution results are not highly sensitive to the specific type of kernel used.

For this reason, we test using a kernel type with varying support sizes. Specifically, we employ a 2D version of Kernel A shown in Fig. \ref{fig:kernels}, a radial kernel that decreases at a constant rate as the distance from the origin increases. Fig. \ref{fig:jump_ind_comp} illustrates the jump indicators defined in \eqref{eq:jumpInd2D}, which are utilized for both edge reclassification and the proposed weight design.
As the kernel support size increases from 3 to 7, oscillations are more effectively flattened, simplifying the distinction between edges and oscillations. However, this comes at the cost of reduced sharpness in the estimated edge locations. Importantly, this reduction in sharpness has minimal impact on the exponent design, as the process ultimately employs patch-wise assignment of exponents. This ensures that an exponent value of 1 is assigned within the vicinity of the edge, maintaining the design's effectiveness.

The VBJS weighting strategy is designed to pinpoint and recover discontinuities in piecewise smooth functions. However, the multifeatured solutions of interest in our work often include large-amplitude oscillations and intricate edges, whose reconstructions can be highly sensitive to the choice of regularization parameters. Consequently, these features exhibit large pointwise variances, resulting in small weights, as shown in Fig. \ref{fig:exponent_weight_ABC}(A3)-(C3). These small weights can disproportionately promote certain edges and may even emphasize oscillatory regions with large amplitudes more strongly than some edges.

In contrast, the proposed weights in Fig. \ref{fig:exponent_weight_ABC}(A4)-(C4) reflect the overall variation of the solution. By focusing on balancing the recovery of multiple features, the proposed weights maintain a relatively narrow maximum range of $[1, 2]$, achieved by setting $(a,b)$ to $(-1,2)$ in Algorithm \ref{alg:newinhomoreg}. This controlled range ensures equitable emphasis across features, as opposed to the broader maximum range of $[0, 100]$ seen with VBJS weights. The wider range of VBJS weights tends to favor the promotion of specific features, particularly those with significantly small weights, potentially at the expense of balanced multifeature recovery.

\begin{figure}[bt!]
	\makebox[\linewidth]{
		\includegraphics[width=\textwidth]{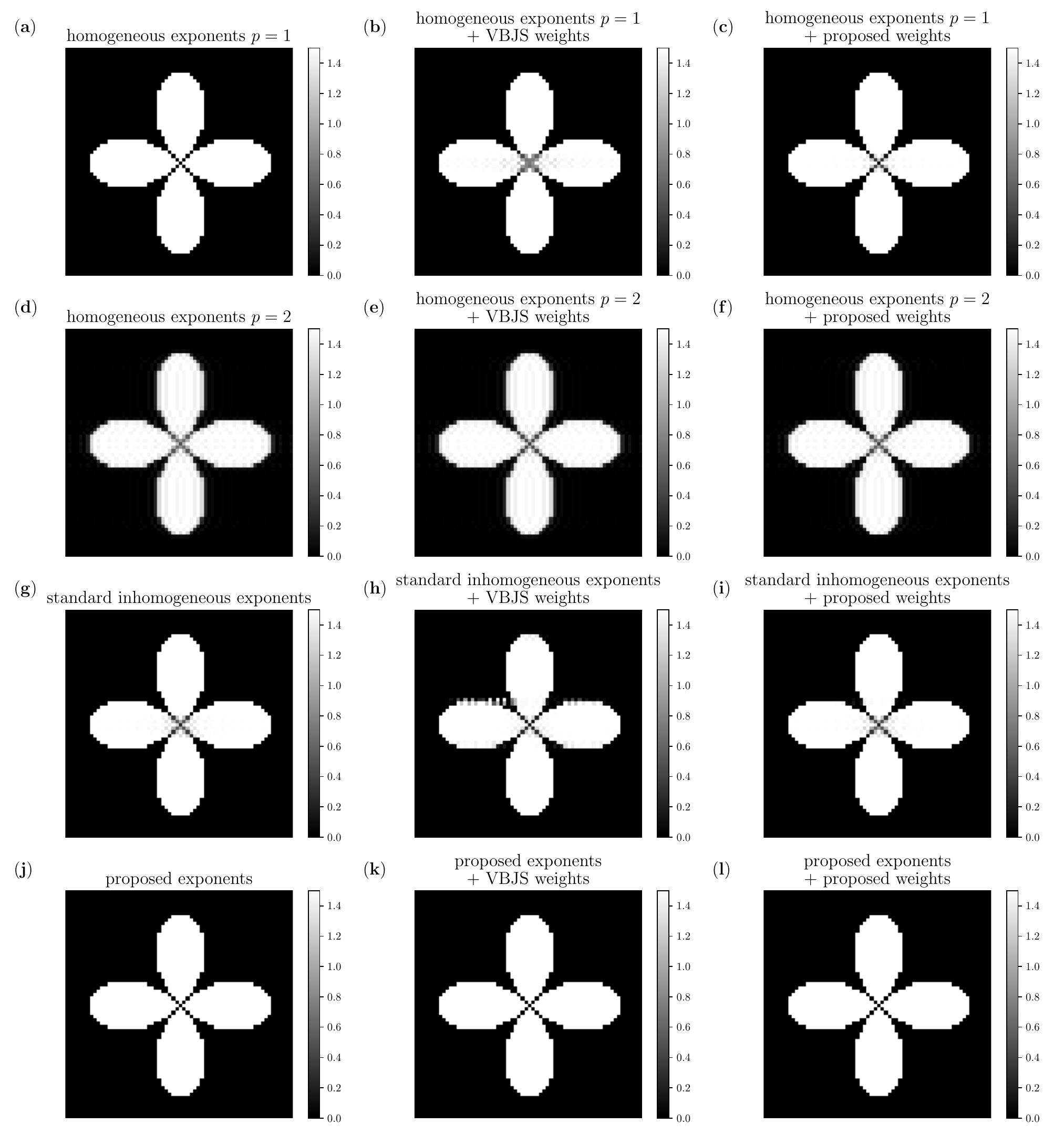}}
	\caption{Reconstructions of Image A in Fig \ref{fig:true_2d_synthetic}(a)
		\label{fig:reconstruction_A}}
\end{figure}

\begin{figure}[bt!]
	\makebox[\linewidth]{
		\includegraphics[width=\textwidth]{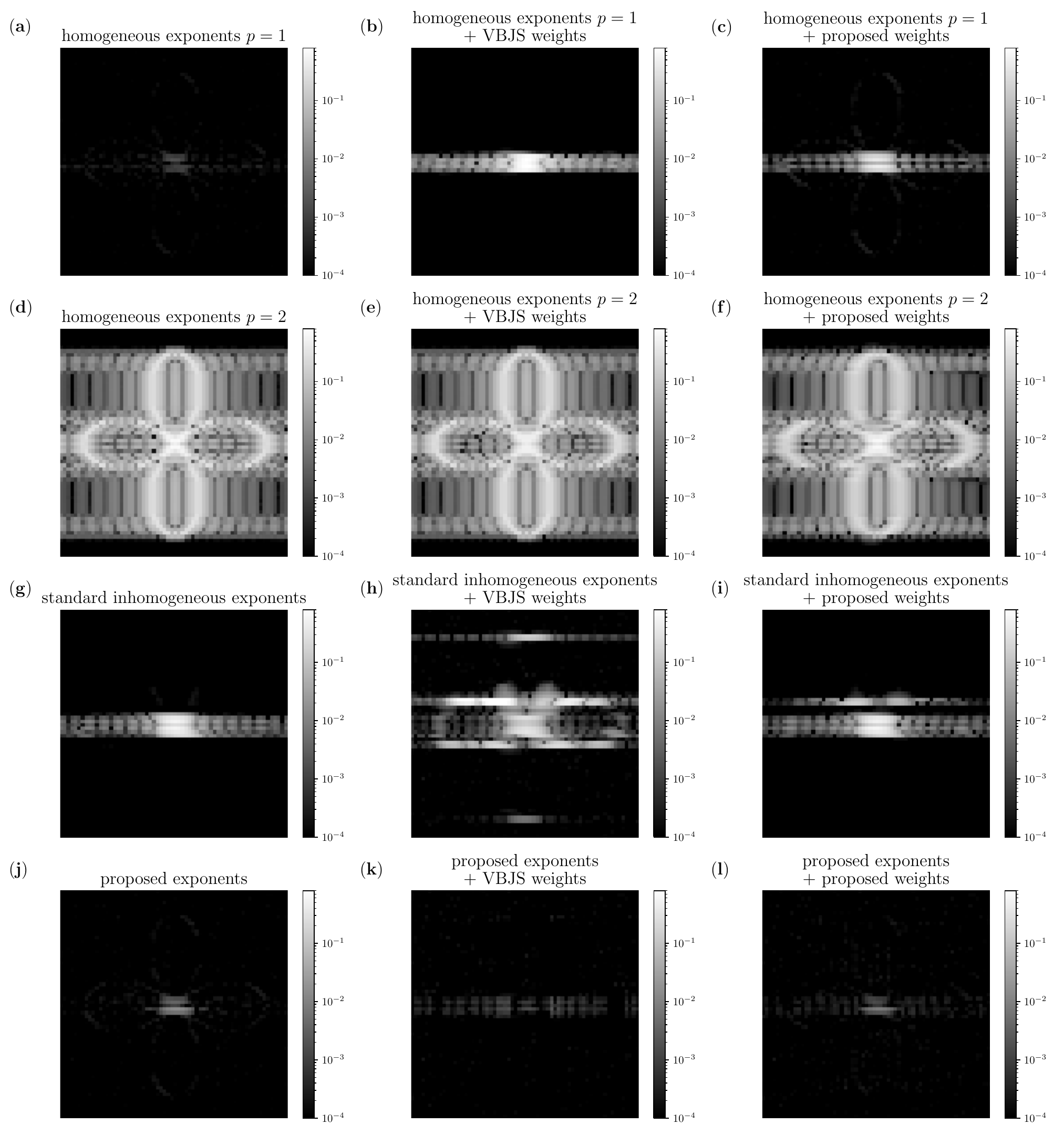}}
	\caption{Reconstruction pointwise errors of Image A in Fig \ref{fig:true_2d_synthetic}(a)
		\label{fig:ptwise_err_A}}
\end{figure}

\begin{figure}[tb!]
	\makebox[\linewidth]{
		\includegraphics[width=\textwidth]{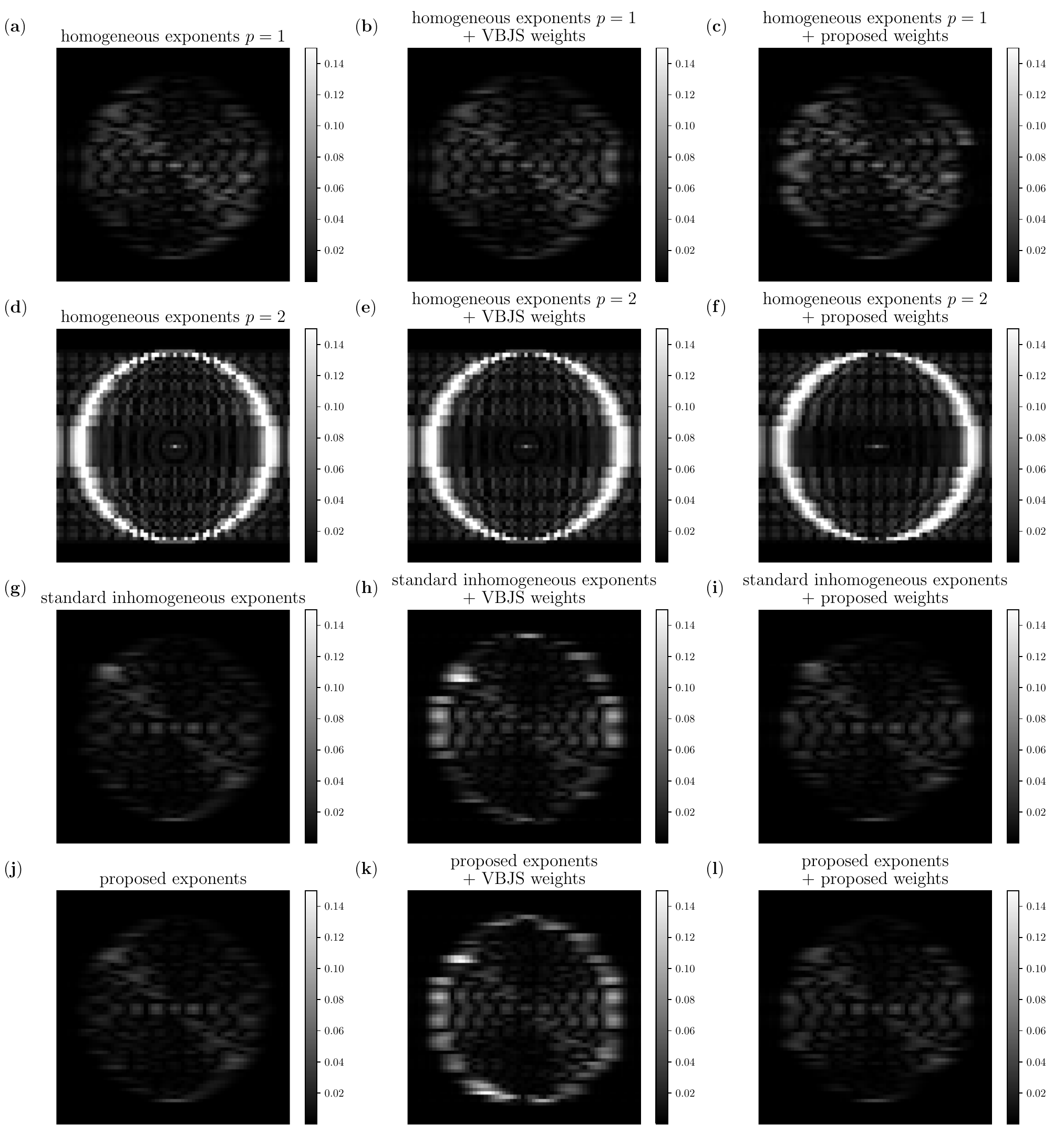}}
	\caption{Reconstruction pointwise errors of Image B in Fig \ref{fig:true_2d_synthetic}(b)
		\label{fig:ptwise_err_B}}
\end{figure}

\begin{figure}[tb!]
	\makebox[\linewidth]{
		\includegraphics[width=\textwidth]{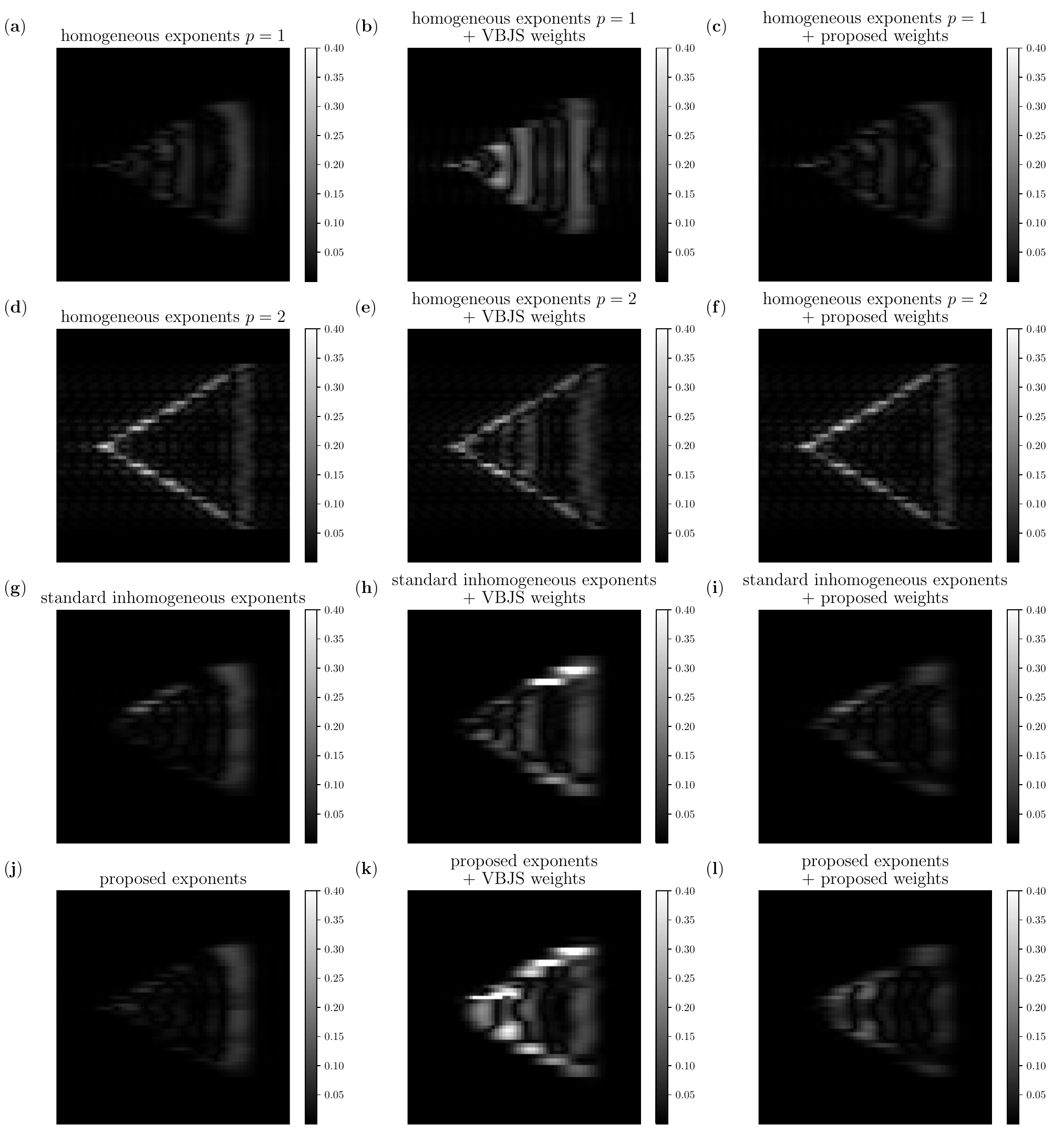}}
	\caption{Reconstruction pointwise errors of Image C in Fig \ref{fig:true_2d_synthetic}(c)
		\label{fig:ptwise_err_C}}
\end{figure}

We evaluate various recovery methods using homogeneous $p=1$ and $p=2$ exponents, as well as standard and proposed inhomogeneous exponents. Additionally, we assess recoveries using these exponents with and without the proposed and alternative weighting strategies. For all tests with Images A-C, we set the parameters in Algorithm \ref{alg:newinhomoreg} to $(l, \tau, a, b) = (3, 0.35, -1, 2)$ for all tests with Images A-C.
As shown in Table \ref{tabl:recon_err2D_imageA_noiseless}, the homogeneous $p=1$ exponent without weights produces the lowest reconstruction errors, as expected for the piecewise constant Image A. While the reconstruction errors using the proposed inhomogeneous exponents are higher than those for $p=1$, they remain significantly lower than those for the standard inhomogeneous exponents due to improved feature classification.
 
Fig. \ref{fig:reconstruction_A} highlights the visual quality of reconstructions. Reconstructions from the unweighted homogeneous $p=1$ regularization (a) and the unweighted and weighted proposed inhomogeneous regularizations (j)-(l) appear similarly accurate. However, their differences become apparent in the pointwise errors shown in Fig. \ref{fig:ptwise_err_A}. The homogeneous $p=1$ regularization (a) exhibits small pointwise errors overall, while the proposed inhomogeneous exponents introduce slightly larger errors in the central region, corresponding to the crossing edges. In contrast, the standard inhomogeneous exponents result in significantly larger pointwise errors in this area due to the misclassification of features and incorrect exponent assignments.
The homogeneous $p=2$ regularization consistently produces the highest reconstruction errors, reflecting its limitations in recovering discontinuities. Therefore, for the remainder of this section, we focus on comparing the performance of the homogeneous $p=1$ and inhomogeneous regularizations.

\begin{table}[tb!]
	\centering
	\footnotesize
	\makebox[\linewidth]{
		\begin{tabular}{|c|c|c|c|}
			\hline
			& relative $\ell_1$ error & relative $\ell_2$ error & relative $\ell_{\infty}$ error \\ \hline
			homogeneous $p=1$                                              & 1.4935e-04          & 1.3884e-04          & 1.7343e-03    \\
			homogeneous $p=1$ + VBJS weight                                & 1.8828e-02          & 7.0314e-02          & 6.5497e-01   \\
			homogeneous $p=1$ + proposed weight                		       & 6.1097e-03          & 2.8594e-02          & 4.1901e-01   \\
			homogeneous $p=2$                                 			   & 1.0647e-01          & 1.0338e-01          & 6.3927e-01    \\
			homogeneous $p=2$ + VBJS weight                   			   & 1.0408e-01          & 9.8110e-02          & 4.7631e-01   \\
			homogeneous $p=2$ + proposed weight              			   & 9.4508e-02          & 9.2675e-02          & 6.1766e-01    \\
			standard inhomogeneous exponent                    & 1.2188e-02          & 4.8407e-02          & 5.6039e-01   \\
			standard inhomogeneous exponent + VBJS weight      & 3.6630e-02          & 1.0274e-01          & 1.3955e+00   \\
			standard inhomogeneous exponent + proposed weight  & 1.2861e-02          & 4.7441e-02          & 5.8129e-01   \\
			proposed inhomogeneous exponent                    & 2.4753e-04          & 6.8543e-04          & 1.2812e-02   \\
			proposed inhomogeneous exponent + VBJS weight      & 1.8345e-04          & 1.9149e-04          & 2.1425e-03   \\
			proposed inhomogeneous exponent + proposed weight  & 2.2023e-04          & 4.0116e-04          & 7.7220e-03   \\
			\hline
		\end{tabular}
	}
	\caption{Reconstruction errors from the recoveries of Image A in Fig \ref{fig:true_2d_synthetic}(a) \label{tabl:recon_err2D_imageA_noiseless}}
\end{table}

\begin{table}[htb!]
	\centering
	\footnotesize
	\makebox[\linewidth]{
		\begin{tabular}{|c|c|c|c|}
			\hline
			& relative $\ell_1$ error & relative $\ell_2$ error & relative $\ell_{\infty}$ error \\ \hline
			homogeneous $p=1$                                              & 1.3820e-02          & 1.6335e-02          & 6.7453e-02    \\
			homogeneous $p=1$ + VBJS weight                                & 1.4845e-02          & 1.7156e-02          & 6.3858e-02   \\
			homogeneous $p=1$ + proposed weight                		       & 1.6138e-02          & 1.9470e-02          & 7.4659e-02    \\
			homogeneous $p=2$                                 			   & 6.2331e-02          & 7.6017e-02          & 1.9757e-01    \\
			homogeneous $p=2$ + VBJS weight                   			   & 6.1360e-02          & 7.4379e-02          & 1.9518e-01   \\
			homogeneous $p=2$ + proposed weight              			   & 5.6043e-02          & 6.7983e-02          & 2.1727e-01    \\
			standard inhomogeneous exponent                    & 8.5299e-03          & 1.1045e-02          & 6.6645e-02   \\
			standard inhomogeneous exponent + VBJS weight      & 1.7209e-02          & 2.3335e-02          & 1.5028e-01   \\
			standard inhomogeneous exponent + proposed weight  & 9.2420e-03          & 1.1918e-02          & 5.8196e-02   \\
			proposed inhomogeneous exponent                    & 8.2371e-03          & 1.0472e-02          & 4.5732e-02   \\
			proposed inhomogeneous exponent + VBJS weight      & 2.1964e-02          & 2.8605e-02          & 1.5366e-01   \\
			proposed inhomogeneous exponent + proposed weight  & 9.4015e-03          & 1.1871e-02          & 4.4392e-02   \\
			\hline
		\end{tabular}
	}
	\caption{Reconstruction errors from the recoveries of Image B in Fig \ref{fig:true_2d_synthetic}(b) \label{tabl:recon_err2D_imageB_noiseless}}
\end{table}

\begin{table}[htb!]
	\centering
	\footnotesize
	\makebox[\linewidth]{
		\begin{tabular}{|c|c|c|c|}
			\hline
			& relative $\ell_1$ error & relative $\ell_2$ error & relative $\ell_{\infty}$ error \\ \hline
			homogeneous $p=1$                                  & 4.1385e-02          & 4.3150e-02          & 9.6352e-02    \\
			homogeneous $p=1$ + VBJS weight                    & 8.2988e-02          & 7.8985e-02          & 2.4084e-01   \\
			homogeneous $p=1$ + proposed weight                & 4.0101e-02          & 4.1541e-02          & 9.7942e-02    \\
			homogeneous $p=2$                                  & 8.6208e-02          & 7.9688e-02          & 3.0480e-01    \\
			homogeneous $p=2$ + VBJS weight                    & 8.9890e-02          & 7.7129e-02          & 2.7357e-01   \\
			homogeneous $p=2$ + proposed weight                & 8.3895e-02          & 7.7272e-02          & 2.9536e-01    \\
			standard inhomogeneous exponent                    & 3.5483e-02          & 4.1950e-02          & 1.5781e-01   \\
			standard inhomogeneous exponent + VBJS weight      & 8.3362e-02          & 1.0830e-01          & 4.9645e-01   \\
			standard inhomogeneous exponent + proposed weight  & 3.3010e-02          & 3.7849e-02          & 1.6981e-01   \\
			proposed inhomogeneous exponent                    & 3.2181e-02          & 3.7675e-02          & 1.0302e-01   \\
			proposed inhomogeneous exponent + VBJS weight      & 1.2169e-01          & 1.5952e-01          & 9.6569e-01   \\
			proposed inhomogeneous exponent + proposed weight  & 3.9472e-02          & 4.6394e-02          & 1.5735e-01   \\
			\hline
		\end{tabular}
	}
	\caption{Reconstruction errors from the recoveries of Image C in Fig \ref{fig:true_2d_synthetic}(c)   \label{tabl:recon_err2D_imageC_noiseless}}
\end{table}

The advantage of inhomogeneous exponents in a target recovery can be validated in the test with Images B and C, which include large amplitude oscillation, edge, and smooth region. However, the two images have differences in the shapes of those features. Especially the fluctuating edges of Image C make the feature classification more difficult. Both standard and proposed inhomogeneous exponents reduce reconstruction errors significantly compared to the homogeneous $p=1$ as shown in Tables \ref{tabl:recon_err2D_imageB_noiseless} and \ref{tabl:recon_err2D_imageC_noiseless}. Moreover, the proposed inhomogeneous exponents show smaller reconstruction errors than the standard inhomogeneous exponents and also reduce the pointwise errors in a part of the edges in the northwest direction as shown in Fig \ref{fig:ptwise_err_B}(j) and Fig \ref{fig:ptwise_err_C}(j) by correcting the missing edge.

For the recoveries of Images A-C, both proposed and VBJS weights are not significantly helpful, and they generally worsen the recovery quality. However, due to VBJS's extreme weight distribution, we can observe its pointwise errors in Figs \ref{fig:ptwise_err_A}, \ref{fig:ptwise_err_B}, and \ref{fig:ptwise_err_C} are often larger compared to the unweighted and the proposed weighted regularizations. These results can be improved by adjusting the choice of $\hat{\epsilon}$ and the weight function in \eqref{eq:vbjs}, resulting in the weight range reduction. In the following subsection, we compare reconstructions of Images D-F obtained by shrinking Images A-C to validate the advantage of weights integrated into the inhomogeneous regularization to reconstruct small features.

\subsubsection{Recovery of scaled images with small features}

\begin{figure}[htb!]
	\makebox[\linewidth]{
		\includegraphics[width=\textwidth]{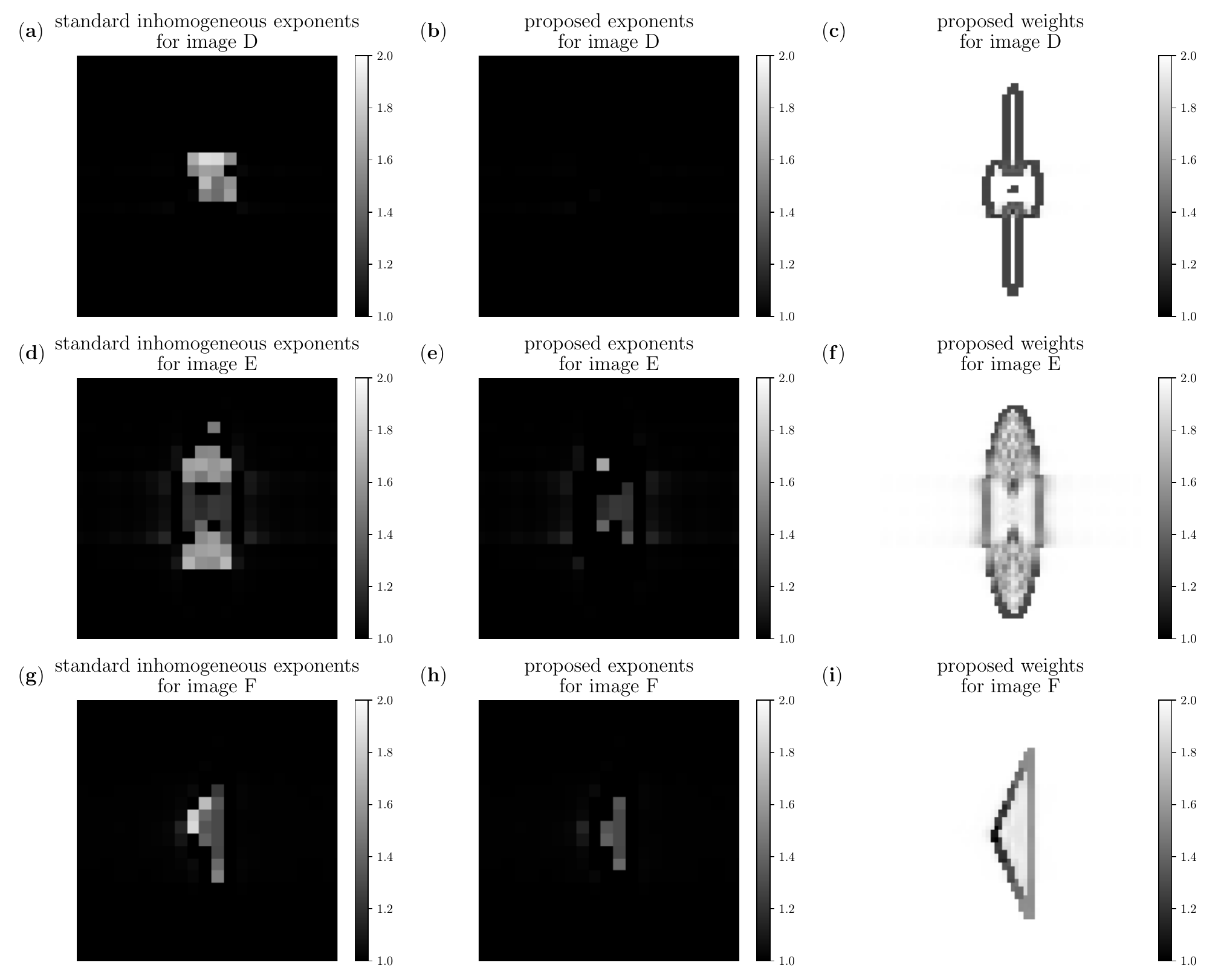}}
	\caption{Standard and proposed inhomogeneous exponents, and proposed weights for Images D-F in Fig \ref{fig:true_2d_synthetic}(d)-(f)
		\label{fig:exponent_weight_DEF}}
\end{figure}

Images D-F in Fig \ref{fig:true_2d_synthetic} have small and narrow features due to scaling, while the fundamental shapes of the features are the same as those in Images A-C, respectively. However, we see in this section that the recovery patterns by various regularizations are not the same, i.e., different regularizations are preferred depending on whether the features are smaller or larger. Particularly, small features require ``weighted" inhomogeneous regularizations to enhance the regularization effect of the small exponent clusters and provide extra information about fine variations within a patch. In contrast, properly designed exponents without weights recover large features well enough.

\begin{figure}[tb!]
	\makebox[\linewidth]{
		\includegraphics[width=\textwidth]{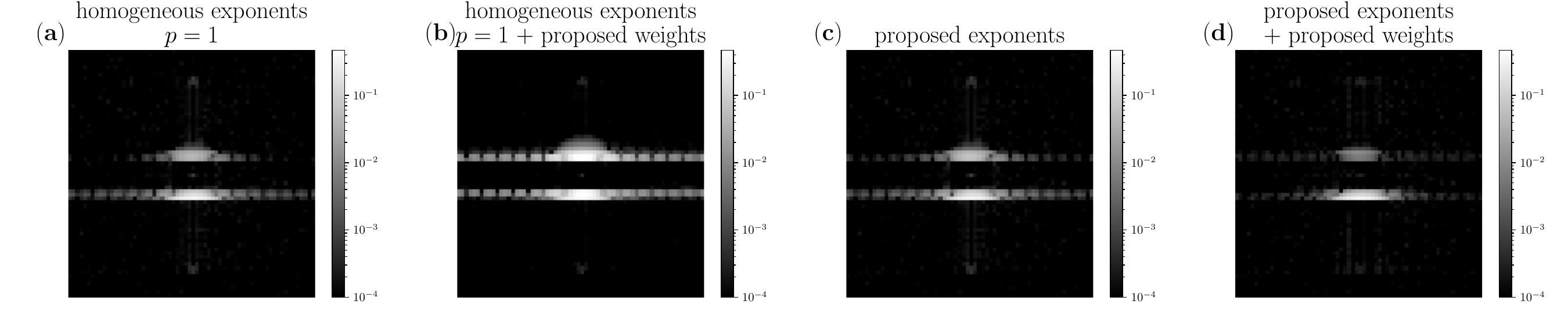}}
	\caption{Reconstruction pointwise errors of Image D in Fig \ref{fig:true_2d_synthetic}(d)
		\label{fig:ptwise_err2d_imageD_noiseless}}
\end{figure}

Image D has isolated discontinuities between each pair of leaves and at the center, which makes consecutive sudden intensity changes. Therefore, even though it is a piecewise constant image as Image A, the homogeneous $p=1$ exponent is insufficient for a good reconstruction. Even with weights in Fig \ref{fig:exponent_weight_DEF}(c) reflecting the intensity variation, homogeneous $p=1$ exponents lead to larger reconstruction and point-wise errors as shown in Table \ref{tabl:recon_err2D_imageD_noiseless} and Fig \ref{fig:ptwise_err2d_imageD_noiseless}(b). However, when we combine the proposed exponents and weights, the reconstruction is improved, achieving the least reconstruction errors. We also observe small point-wise errors, especially near the isolated edges in the center and edges along the narrow leaves pointing north and south, as seen in Fig \ref{fig:ptwise_err2d_imageD_noiseless}(d). This result demonstrates that the weighted inhomogeneous regularization is good at feature recovery when a single patch contains multiple features. The patch-wise exponent assignment lacks information on the variation within a patch, but the point-wise weights complement this information. The proposed inhomogeneous exponents without weights give reconstruction errors larger than the homogeneous $p=1$ even though the exponents in Fig \ref{fig:exponent_weight_DEF}(a) are correctly designed so that the exponents are 1 or close to 1. Therefore, the weights integrated into the inhomogeneous regularization are necessary to make the inhomogeneous exponents effective in a small region.

\begin{figure}[tb!]
	\makebox[\linewidth]{
		\includegraphics[width=\textwidth]{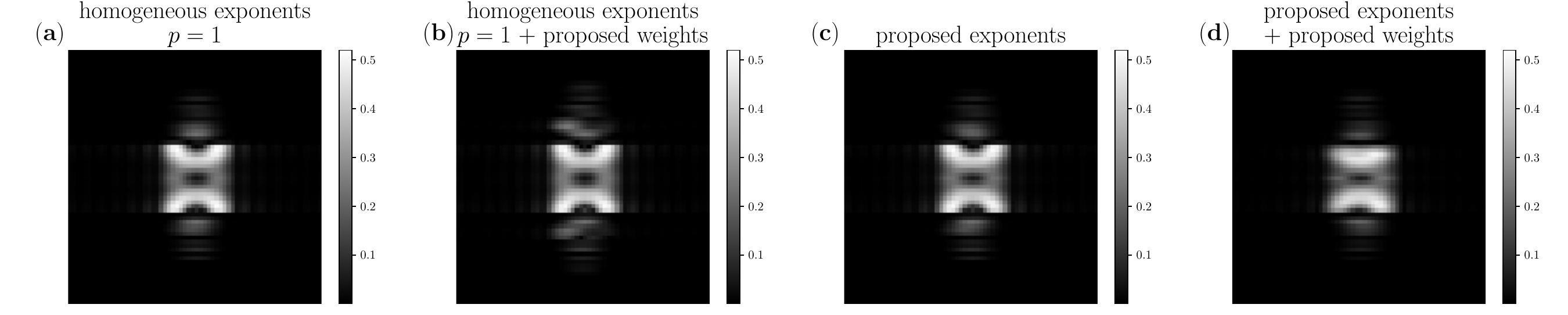}}
	\caption{Reconstruction pointwise errors of Image E in Fig \ref{fig:true_2d_synthetic}(e)
		\label{fig:ptwise_err2d_imageE_noiseless}}
\end{figure}

\begin{figure}[tb!]
	\makebox[\linewidth]{
		\includegraphics[width=\textwidth]{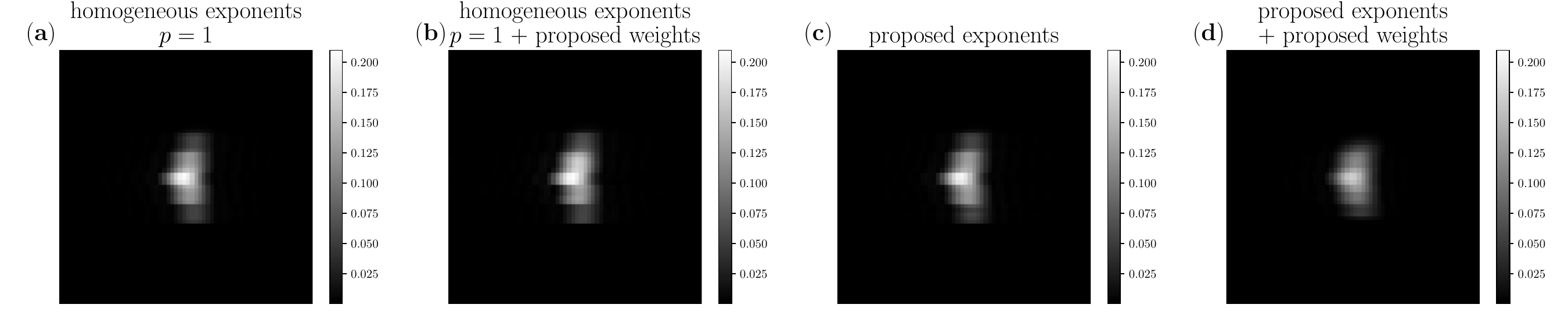}}
	\caption{Reconstruction pointwise errors of Image F in Fig \ref{fig:true_2d_synthetic}(f)
		\label{fig:ptwise_err2d_imageF_noiseless}}
\end{figure}

For Images E and F with multiple narrow features, the proposed exponents in Fig \ref{fig:exponent_weight_DEF}(e) and (h) appropriately assign 1 at the edge and relatively larger values for the oscillation inside while the standard inhomogeneous regularizations in Fig \ref{fig:exponent_weight_DEF}(d) and (g) incorrectly assign relatively large exponents for the regions including both edges and oscillations. The proposed weights in Fig \ref{fig:exponent_weight_DEF}(f) and (i) also reflect the targets' details closely.
Even though Images E and F have multiple features as Images B and C, inhomogeneous exponents without weights do not significantly reduce the reconstruction errors compared to the homogeneous $p=1$ exponents as in Tables \ref{tabl:recon_err2D_imageE_noiseless} and \ref{tabl:recon_err2D_imageF_noiseless}. However, we achieve a larger error reduction after incorporating the weights. Pointwise errors in Figs \ref{fig:ptwise_err2d_imageE_noiseless} and \ref{fig:ptwise_err2d_imageF_noiseless}(d) also show smaller errors along the edge and over the oscillation inside. These results validate that integrating proposed inhomogeneous exponents and weights works well in reconstructing a small area of a single feature type by promoting the corresponding small exponent cluster's influence.

Through the tests with Images A-F in this section, we can validate that the proposed inhomogeneous exponents correctly classify features of various shapes and sizes, and the proposed weights portray the detailed variations in each target image. The inhomogeneous exponent alone performs well in recovering the mixed features of sufficiently large sizes, as in Images B and C. However, a small area of single feature and mixed features within a patch, as in scaled Images D-F, further require weights to be combined into the inhomogeneous regularization to enhance the efficacy of the inhomogeneous exponents.

\begin{table}[tb!]
	\centering
	\footnotesize
	\makebox[\linewidth]{
		\begin{tabular}{|c|c|c|c|}
			\hline
			& relative $\ell_1$ error & relative $\ell_2$ error & relative $\ell_{\infty}$ error \\ \hline
			homogeneous $p=1$      									& 1.4634e-02          & 3.8056e-02          & 3.2641e-01    \\ 			
			homogeneous $p=1$  + proposed weights      				& 5.7877e-02          & 1.0501e-01          & 6.2030e-01    \\ 
			standard inhomogeneous exponent                   		& 5.9792e-02          & 1.0002e-01          & 5.2582e-01    \\
			standard inhomogeneous exponent + proposed weight 		& 4.8655e-02          & 8.7655e-02          & 4.6886e-01    \\
			proposed inhomogeneous exponents                		& 1.6075e-02          & 3.9323e-02          & 3.2963e-01    \\  
			proposed inhomogeneous exponents + proposed weights 	& 1.1660e-02          & 3.6052e-02          & 3.2507e-01    \\ 
			\hline
		\end{tabular}
	}
	\caption{Reconstruction errors from the recoveries of Image D in Fig \ref{fig:true_2d_synthetic}(d)   \label{tabl:recon_err2D_imageD_noiseless}}
\end{table}

\begin{table}[tb!]
	\centering
	\footnotesize
	\makebox[\linewidth]{
		\begin{tabular}{|c|c|c|c|}
			\hline
			& relative $\ell_1$ error & relative $\ell_2$ error & relative $\ell_{\infty}$ error \\ \hline
			homogeneous $p=1$      									& 2.1249e-01          & 2.2650e-01          & 5.0407e-01    \\ 			
			homogeneous $p=1$  + proposed weights      				& 2.3195e-01          & 2.3078e-01          & 5.0662e-01     \\ 
			standard inhomogeneous exponent                         & 2.1331e-01          & 2.2163e-01          & 4.9990e-01   \\
			standard inhomogeneous exponent + proposed weight       & 2.0523e-01          & 2.2203e-01          & 5.0594e-01    \\
			proposed inhomogeneous exponents                		& 2.1183e-01          & 2.2539e-01          & 5.0533e-01    \\  
			proposed inhomogeneous exponents + proposed weights 	& 1.8687e-01          & 2.0941e-01          & 4.9156e-01    \\ 
			\hline
		\end{tabular}
	}
	\caption{Reconstruction errors from the recoveries of Image E in Fig \ref{fig:true_2d_synthetic}(e)    \label{tabl:recon_err2D_imageE_noiseless}}
\end{table}

\begin{table}[tb!]
	\centering
	\footnotesize
	\makebox[\linewidth]{
		\begin{tabular}{|c|c|c|c|}
			\hline
			& relative $\ell_1$ error & relative $\ell_2$ error & relative $\ell_{\infty}$ error \\ \hline
			homogeneous $p=1$                                  & 9.4195e-02          & 9.0923e-02          & 1.9379e-01     \\
			homogeneous $p=1$ + proposed weight                & 1.0366e-01          & 9.8911e-02          & 2.0753e-01    \\
			standard inhomogeneous exponent                    & 9.4949e-02          & 9.3681e-02          & 2.3550e-01   \\
			standard inhomogeneous exponent + proposed weight  & 8.0469e-02          & 8.1907e-02          & 2.1293e-01    \\
			proposed inhomogeneous exponent                    & 9.4652e-02          & 9.0514e-02          & 2.0046e-01   \\
			proposed inhomogeneous exponent + proposed weight  & 7.4409e-02          & 7.6126e-02          & 1.7102e-01    \\
			\hline
		\end{tabular}
	}
	\caption{Reconstruction errors from the recoveries of Image F in Fig \ref{fig:true_2d_synthetic}(f)    \label{tabl:recon_err2D_imageF_noiseless}}
\end{table}

\subsection{Real sea ice image recovery}\label{subsec:exp:seaice}

\begin{figure}[tb!]
	\makebox[\linewidth]{
		\includegraphics[width=\textwidth]{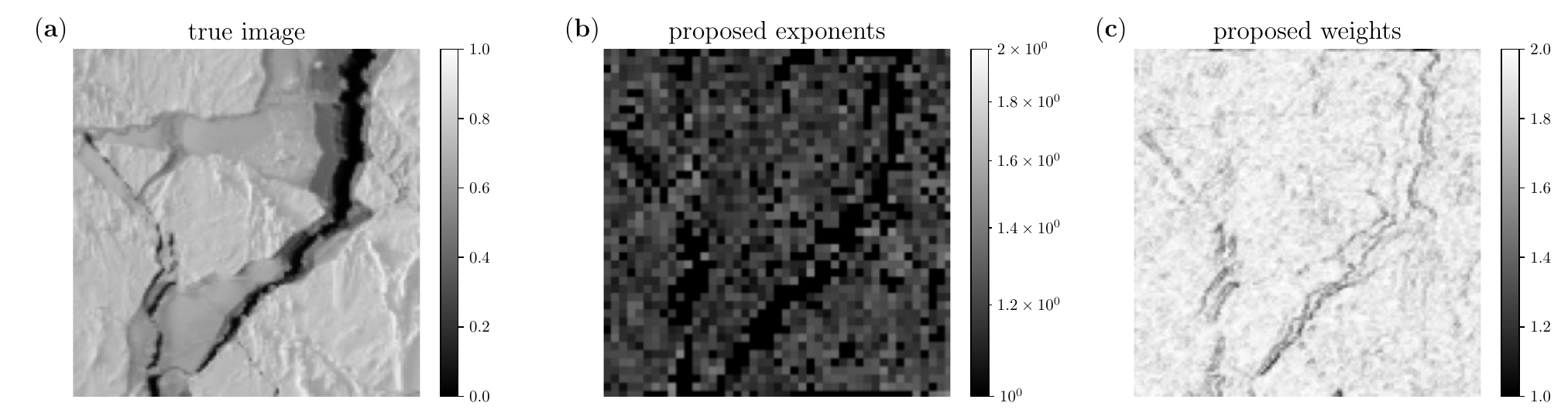}}
	\caption{2D sea ice image, and proposed exponents and weights }
	\label{fig:true_exponent_weight_comparison_2d_seaice_noiseless}
\end{figure}

\begin{table}[tb!]
	\centering
	\footnotesize
	\makebox[\linewidth]{
		\begin{tabular}{|c|c|c|c|}
			\hline
			& relative $\ell_1$ error & relative $\ell_2$ error & relative $\ell_{\infty}$ error \\ \hline
			homogeneous $p=1$      								& 8.8431e-02          & 1.1286e-01          & 3.4506e-01   \\ 			
			homogeneous $p=1$  + proposed weights      			& 9.1763e-02          & 1.1832e-01          & 3.6380e-01   \\ 
			proposed inhomogeneous exponents                	& 9.5969e-02          & 1.3066e-01          & 4.0536e-01   \\  
			proposed inhomogeneous exponents + proposed weights & 8.5859e-02          & 1.0636e-01          & 3.1832e-01   \\ 
			\hline
		\end{tabular}
	}
	\caption{Reconstruction errors from the recoveries of real sea ice image in Fig \ref{fig:true_exponent_weight_comparison_2d_seaice_noiseless}(a) using noisy measurement  \label{tabl:recon_err2D_seaice_noiseless}}
\end{table}

In addition to testing with 2D synthetic data, we apply the proposed regularization to a real sea ice image obtained from the NASA Earth Observatory \cite{nasa}. Direct measurements of sea ice characteristics, such as optical and in-situ data, are expensive to obtain. Therefore, noisy indirect data from satellites, such as ICESat-2 and CryoSat-2, are more commonly used for sea ice observation. As discussed in Section \ref{sec:intro}, the surface of sea ice exhibits various features and textures resulting from the growth and deformation of the ice and snow, which vary in thickness. Sea ice researchers rely on image intensities from remote sensing to estimate important parameters such as sea ice thickness, snow density, and freeboard. These parameters are crucial for understanding the age and dynamics of sea ice.

In the real sea ice image shown in Fig. \ref{fig:true_exponent_weight_comparison_2d_seaice_noiseless}(a), we observe various features and textures, including cracks and elevation variations on the sea ice surface. These images provide crucial information for modeling sea ice dynamics. Therefore, it is essential to obtain high-quality images from remote sensing data that preserve these features. For the experiment, we select a portion of the image and downsample it to a $128\times128$ resolution. Noisy measurements are generated by selecting every third Fourier coefficient in the $x$-direction and adding Gaussian noise from $\mathcal{N}(\mathbf{0}, \sigma^2 I)$. The noise standard deviation, $\sigma$, is set to 3.45, resulting in 25 signal-to-noise ratio (SNR) in decibels (${\rm SNR_{db}}$). For this experiment, we use the parameters $(l, \tau, a, b) = (3, 0.35, -1, 2)$ in Algorithm \ref{alg:newinhomoreg}.

\begin{figure}[tb!]
	\makebox[\linewidth]{
		\includegraphics[width=\textwidth]{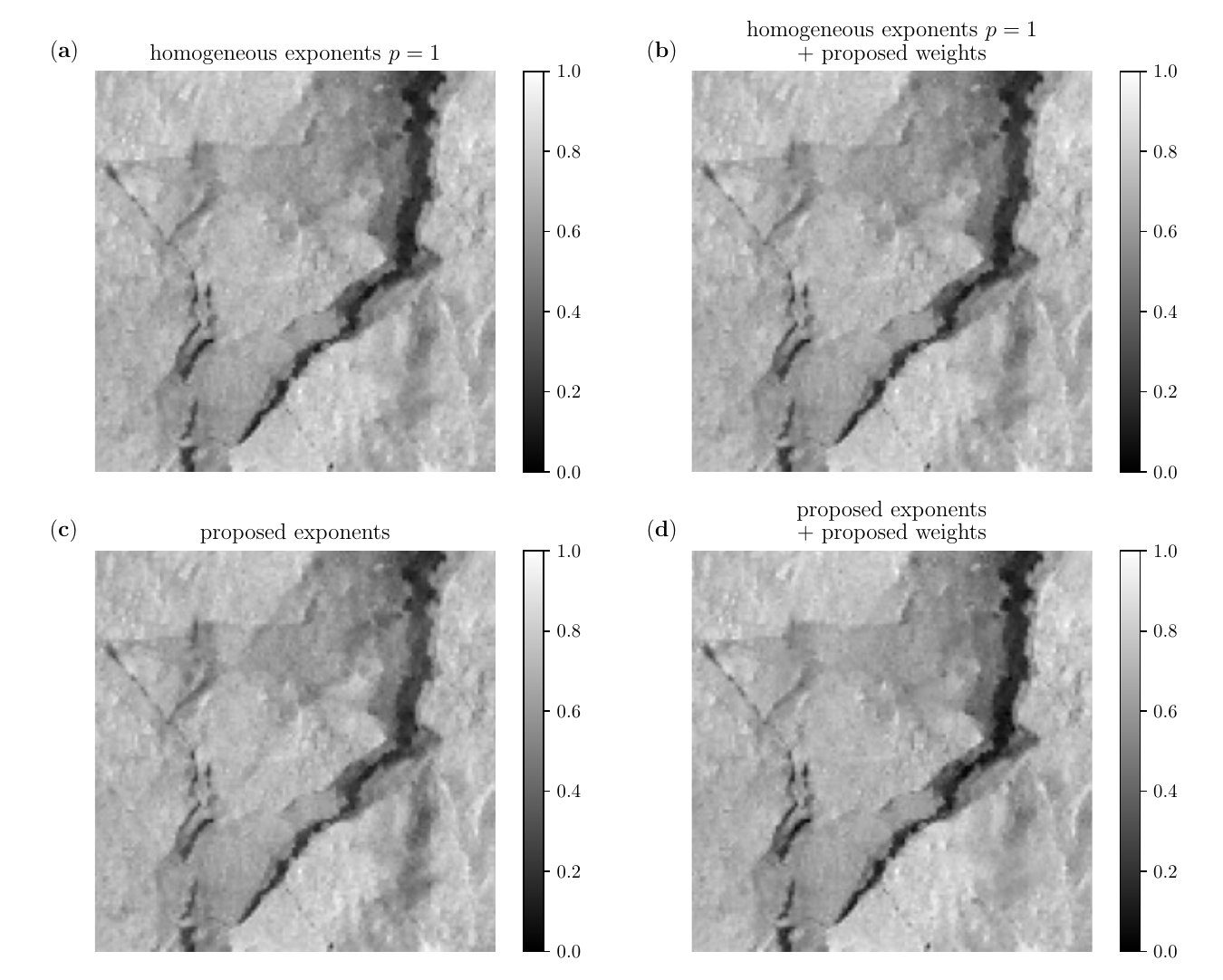}}
	\caption{Reconstructions of the real sea ice image in Fig \ref{fig:true_exponent_weight_comparison_2d_seaice_noiseless}(a) using noisy measurement
		\label{fig:recon1d_noisy_2d_seaice_noiseless}}
\end{figure}

\begin{figure}[tb!]
	\makebox[\linewidth]{
		\includegraphics[width=\textwidth]{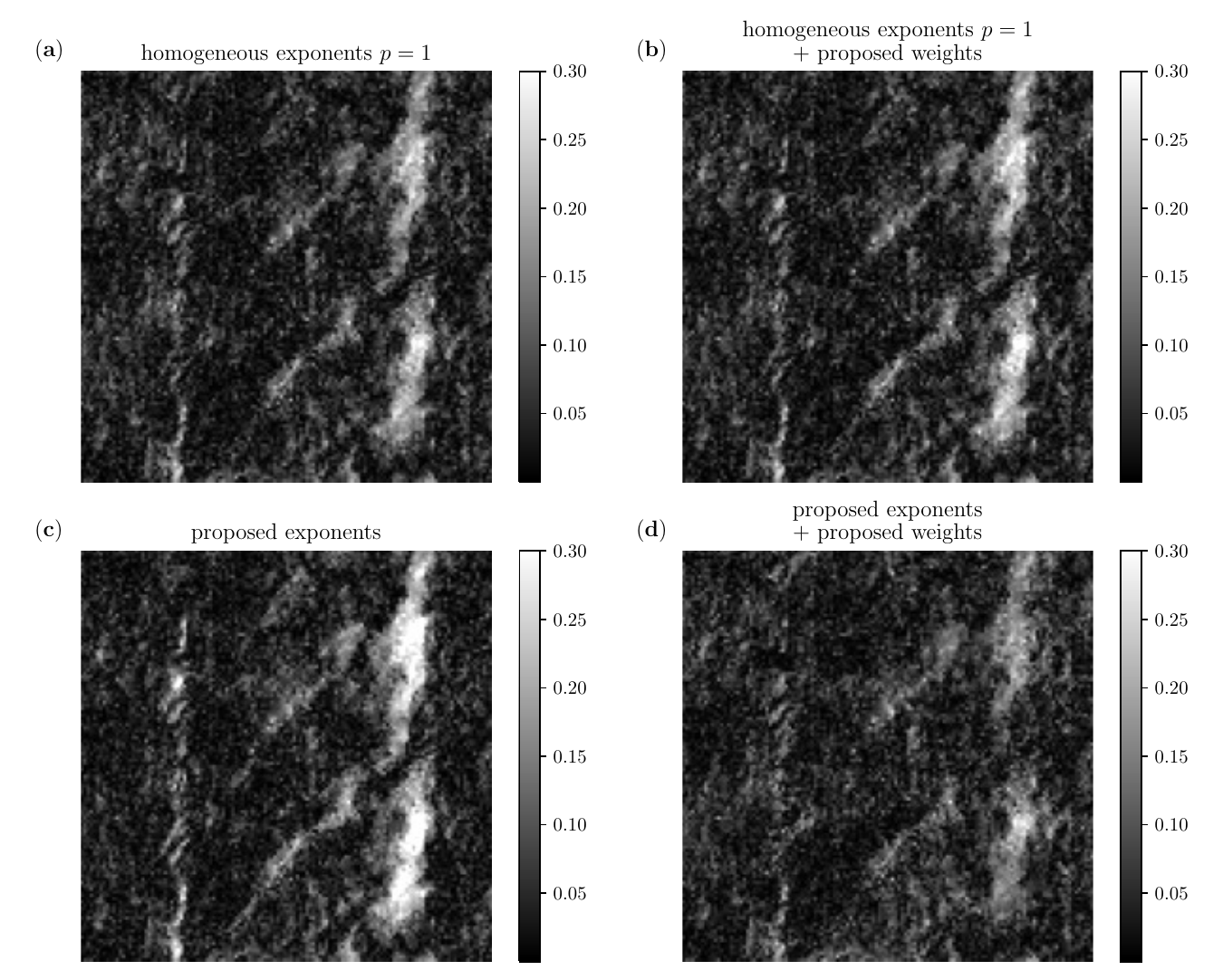}}
	\caption{Recontruction pointwise errors of the real sea ice image in Fig \ref{fig:true_exponent_weight_comparison_2d_seaice_noiseless}(a) using noisy measurement
		\label{fig:ptwise_err2d_seaice_noiseless}}
\end{figure}

Compared to synthetic images, the real sea ice image features multiple, more complex shapes, which complicates the estimation of variations. This complexity is evident in the proposed weight distribution shown in Fig. \ref{fig:true_exponent_weight_comparison_2d_seaice_noiseless}(c), which roughly outlines the edges from cracks and the intensity variations on the surface. We observe that using the homogeneous $p=1$ regularization with the proposed weights leads to a deterioration in recovery quality compared to the unweighted $p=1$ regularization, in terms of both reconstruction and pointwise errors, as seen in Figs \ref{fig:recon1d_noisy_2d_seaice_noiseless}(a) and \ref{fig:ptwise_err2d_seaice_noiseless}(b), respectively.

Table \ref{tabl:recon_err2D_seaice_noiseless} further confirms that the reconstruction errors for $p=1$ with the proposed weights are larger than those for the homogeneous $p=1$ regularization. However, the homogeneous $p=1$ regularization still shows unsatisfactory recovery, as it produces unwanted artifacts, such as a false crack and incorrect intensity values in certain regions, as shown in Fig. \ref{fig:recon1d_noisy_2d_seaice_noiseless}(a). This is verified by the pointwise error in Fig. \ref{fig:ptwise_err2d_seaice_noiseless}(a), which highlights significant errors along the vertical shades on the left and right of the image.

Due to the complex mixture of features in the real sea ice image, the proposed inhomogeneous exponent in Fig. \ref{fig:true_exponent_weight_comparison_2d_seaice_noiseless}(b) displays a more intricate distribution compared to the exponents for the synthetic 2D images in Section \ref{subsec:synthetic2d}. Despite this, it successfully captures important features. For example, it highlights black strips in two directions that correspond to the cracks in the true image. However, in the proposed exponents for the real sea ice image, we observe a small or thin exponent cluster that is isolated and surrounded by clusters of different exponents. Although the inhomogeneous exponent is designed appropriately, both reconstruction errors and pointwise errors are still larger than those for the homogeneous $p=1$ regularization without the weights, as shown in Table \ref{tabl:recon_err2D_seaice_noiseless} and Fig. \ref{fig:ptwise_err2d_seaice_noiseless}(c).

Therefore, the unweighted inhomogeneous regularization fails to accurately reconstruct the various features in the sea ice image. However, by combining the proposed weight and inhomogeneous exponent, we mitigate the drawbacks of both the homogeneous and unweighted inhomogeneous regularizations. This is evident in the reconstruction shown in Fig. \ref{fig:recon1d_noisy_2d_seaice_noiseless}(d), which accurately reflects the true cracks, intensity variations, and surface texture without introducing fake cracks or features. Additionally, the pointwise error in Fig. \ref{fig:ptwise_err2d_seaice_noiseless}(d) is relatively small across the entire domain, indicating that the various characteristics of the sea ice are well recovered.

\section{Discussions and conclusions}\label{sec:conclusion}
This work introduces a weighted inhomogeneous regularization framework for solving inverse problems with multifeatured solutions. Our method leverages both pixel-wise and patch-wise information for feature identification, addressing the limitations of prior inhomogeneous regularizations that rely solely on patch-based statistics. By incorporating pixel-wise information, specifically through the (directional) derivatives of smoothed reconstructions, we mitigate feature misclassification in closely intertwined regions. The smoothing kernel enhances edge separation by suppressing variations in other features, leading to a more accurate inhomogeneous exponent design. To address insufficient regularization in regions with small clusters of constant exponents, we introduce weights that balance feature recovery and improve the representation of variations within a target solution.

We developed a solver based on a modified Alternating Direction Method of Multipliers (ADMM) to handle the weighted inhomogeneous regularization, enabling efficient optimization by separating the fidelity and regularization terms. Numerical experiments with synthetic and real sea ice images demonstrate the robustness of our approach in recovering multifeatured solutions, showing significant improvements in feature classification and reconstruction accuracy.

While our focus has been on partial Fourier measurements due to their relevance in remote sensing applications for sea ice research, other data types, such as in-situ and optical measurements, offer complementary information despite their limited availability due to high costs and collection challenges \cite{insitu, optical}. Future work will explore data fusion strategies to integrate disparate data types effectively, addressing the potential ill-posedness introduced by combining multiple measurement operators \cite{datafusion}.

Our regularization design relies on tuning parameters, such as kernel support and threshold sizes, which were hand-tuned in this study. Future research will involve a rigorous sensitivity analysis to optimize these parameters, particularly for datasets with extreme scale differences among features. Additionally, while this study restricted the regularization exponent range to $[1, 2]$ for convexity, preliminary tests suggest that alternative exponents could yield better results for specific signal types. Investigating such extensions is a promising direction for future research.

Weighted inhomogeneous regularization also has potential applications beyond the scope of this study. For instance, variational data assimilation methods like 3D-Var \cite{3DVar}, which typically use $\ell_2$ regularization with Gaussian priors, could benefit from the flexibility of inhomogeneous regularization to incorporate more complex, non-Gaussian priors. This approach could enhance reconstruction accuracy for a broader range of Bayesian inference problems. In particular, sea ice dynamics evolve over various time scales, from seasonal to climatological changes, necessitating both spatial and temporal regularization for accurate predictions. Weighted $\ell_2$ regularization has been applied dynamically in data assimilation to capture crack propagation in sea ice \cite{LI2024113396}, but it shows limitations in recovering other features like ridges. Combining weighted inhomogeneous regularization with data assimilation could provide a more comprehensive approach, capturing the full range of sea ice variations with higher accuracy. In particular, integrating dynamically designed inhomogeneous regularization could address the non-Gaussian nature of evolving multifeatured sea ice formations.

\subsection*{Acknowledgement}
JH and YL are supported by ONR MURI N00014-20-1-2595.

\bibliographystyle{abbrv}
\bibliography{inhomoregv2}

\begin{thebibliography}{10}

\bibitem{ICESAT2}
W.~Abdalati, H.~J. Zwally, R.~Bindschadler, B.~Csatho, S.~L. Farrell, H.~A.
  Fricker, D.~Harding, R.~Kwok, M.~Lefsky, T.~Markus, A.~Marshak, T.~Neumann,
  S.~Palm, B.~Schutz, B.~Smith, J.~Spinhirne, and C.~Webb.
\newblock The {ICES}at-2 laser altimetry mission.
\newblock {\em Proceedings of the IEEE}, 98(5):735--751, 2010.

\bibitem{VBJS}
B.~Adcock, A.~Gelb, G.~Song, and Y.~Sui.
\newblock Joint sparse recovery based on variances.
\newblock {\em SIAM Journal on Scientific Computing}, 41(1):A246--A268, 2019.

\bibitem{archibald2005polynomial}
R.~Archibald, A.~Gelb, and J.~Yoon.
\newblock Polynomial fitting for edge detection in irregularly sampled signals
  and images.
\newblock {\em SIAM journal on numerical analysis}, 43(1):259--279, 2005.

\bibitem{inhomo_denoising}
P.~Blomgren and T.~Chan.
\newblock Extensions to total variation denoising.
\newblock {\em Proceedings of SPIE - The International Society for Optical
  Engineering}, 3162, 03 1999.

\bibitem{ADMM}
S.~Boyd, N.~Parikh, and E.~Chu.
\newblock {\em Distributed optimization and statistical learning via the
  alternating direction method of multipliers}.
\newblock Now Publishers Inc, 2011.

\bibitem{CS2006}
E.~J. Candes, J.~K. Romberg, and T.~Tao.
\newblock Stable signal recovery from incomplete and inaccurate measurements.
\newblock {\em Communications on Pure and Applied Mathematics: A Journal Issued
  by the Courant Institute of Mathematical Sciences}, 59(8):1207--1223, 2006.

\bibitem{chandrupatla}
T.~R. Chandrupatla.
\newblock A new hybrid quadratic/bisection algorithm for finding the zero of a
  nonlinear function without using derivatives.
\newblock {\em Advances in Engineering Software}, 28(3):145--149, 1997.

\bibitem{inhomo_denoising2010}
Q.~Chen, P.~Montesinos, Q.~S. Sun, P.~A. Heng, and D.~S. Xia.
\newblock Adaptive total variation denoising based on difference curvature.
\newblock {\em Image and Vision Computing}, 28(3):298--306, 2010.

\bibitem{3DVar}
G.~Evensen, F.~C. Vossepoel, and P.~J. van Leeuwen.
\newblock {\em Kalman filters and 3{DV}ar}, pages 63--71.
\newblock Springer International Publishing, Cham, 2022.

\bibitem{CS2013}
S.~Foucart and H.~Rauhut.
\newblock {\em A mathematical introduction to compressive sensing}.
\newblock Birkh\"{a}user Basel, 2013.

\bibitem{seaicethickness}
C.~Haas, J.~Lobach, S.~Hendricks, L.~Rabenstein, and A.~Pfaffling.
\newblock Helicopter-borne measurements of sea ice thickness, using a small and
  lightweight, digital {EM} system.
\newblock {\em Journal of Applied Geophysics}, 67(3):234--241, 2009.
\newblock Airborne Geophysics.

\bibitem{inhomoregv1}
J.~Han and Y.~Lee.
\newblock Inhomogeneous regularization with limited and indirect data.
\newblock {\em Journal of Computational and Applied Mathematics}, 428:115193,
  2023.

\bibitem{cryosat2}
S.~Hendricks, R.~Ricker, and S.~Paul.
\newblock Product user guide \& algorithm specification: {AWI} {C}ryo{S}at-2
  sea ice thickness (version 2.4).
\newblock 2021.

\bibitem{SAR}
C.~V. Jakowatz, D.~E. Wahl, P.~H. Eichel, D.~C. Ghiglia, and P.~A. Thompson.
\newblock {\em Spotlight-mode synthetic aperture radar: a signal processing
  approach: a signal processing approach}.
\newblock Springer Science \& Business Media, 2012.

\bibitem{cryosat2resolution}
N.~T. Kurtz, N.~Galin, and M.~Studinger.
\newblock An improved {C}ryo{S}at-2 sea ice freeboard retrieval algorithm
  through the use of waveform fitting.
\newblock {\em The Cryosphere}, 8(4):1217--1237, 2014.

\bibitem{datafusion}
C.~König, T.~König, S.~Singha, A.~Frost, and S.~Jacobsen.
\newblock Combined use of space borne optical and {SAR} data to improve
  knowledge about sea ice for shipping.
\newblock {\em Remote Sensing}, 13(23), 2021.

\bibitem{LI2024113396}
T.~Li, A.~Gelb, and Y.~Lee.
\newblock A structurally informed data assimilation approach for nonlinear
  partial differential equations.
\newblock {\em Journal of Computational Physics}, 519:113396, 2024.

\bibitem{nasa}
{\relax NASA Earth Observatory}.
\newblock Flying over arctic sea ice, 2011.

\bibitem{insitu}
C.~Perron, C.~Katlein, S.~Lambert-Girard, E.~Leymarie, L.-P. Guinard,
  P.~Marquet, and M.~Babin.
\newblock Development of a diffuse reflectance probe for in situ measurement of
  inherent optical properties in sea ice.
\newblock {\em The Cryosphere}, 15(9):4483--4500, 2021.

\bibitem{optical}
C.~Pohl, L.~Istomina, S.~Tietsche, E.~J\"akel, J.~Stapf, G.~Spreen, and
  G.~Heygster.
\newblock Broadband albedo of {A}rctic sea ice from {MERIS} optical data.
\newblock {\em The Cryosphere}, 14(1):165--182, 2020.

\bibitem{seaice}
D.~A. Rothrock.
\newblock The energetics of the plastic deformation of pack ice by ridging.
\newblock {\em Journal of Geophysical Research (1896-1977)}, 80(33):4514--4519,
  1975.

\bibitem{TV1992}
L.~I. Rudin, S.~Osher, and E.~Fatemi.
\newblock Nonlinear total variation based noise removal algorithms.
\newblock {\em Physica D: nonlinear phenomena}, 60(1-4):259--268, 1992.

\bibitem{SARbook}
M.~Soumekh.
\newblock {\em Synthetic aperture radar signal processing}, volume~7.
\newblock New York: Wiley, 1999.

\bibitem{groupLasso}
M.~Yuan and Y.~Lin.
\newblock Model selection and estimation in regression with grouped variables.
\newblock {\em Journal of the Royal Statistical Society: Series B (Statistical
  Methodology)}, 68(1):49--67, 2006.

\bibitem{seaicemodeling}
J.~Zhang and D.~A. Rothrock.
\newblock Modeling global sea ice with a thickness and enthalpy distribution
  model in generalized curvilinear coordinates.
\newblock {\em Monthly Weather Review}, 131(5):845 -- 861, 2003.

\end{thebibliography}

\end{document}